\documentclass[final]{article}
\usepackage[utf8]{inputenc}

\usepackage{amsmath,amssymb,amsthm}
\usepackage{mathrsfs}
\usepackage{float}
\usepackage{hyperref}
\usepackage{xcolor}
\usepackage{tikz}
\usepackage{subfigure}
\usepackage[numbers,sort&compress]{natbib}
\numberwithin{equation}{section}
\usetikzlibrary{arrows,shapes,chains}
\usepackage{graphicx,color,xcolor}
\usepackage{caption}
\usepackage{subfigure}

\usepackage{epstopdf}
\usepackage{algorithm}
\usepackage{algpseudocode}

\title{A Deep Learning Galerkin Method for the Closed-Loop Geothermal System}
\author{Wen Zhang, {Jian Li \thanks{email: jiaanli@mail.xjtu.edu.cn, jiaaanli@gmail.com, ZW15290229577@163.com. Supported in part by NSF of China (No. 11771259), special support program to develop innovative talents in the region of Shaanxi province, innovation team on computationally efficient numerical methods based on new energy problems in Shaanxi province, and innovative team project of Shaanxi Provincial Department of Education(No. 21JP013). Department of Mathematics, Shaanxi University of Science and Technology, Xi'an 710021, China.}}
}
\date{}
\begin{document}

\maketitle
\noindent{\bf Abstract:~~} There has been an arising trend of adopting deep learning methods to study partial differential equations (PDEs). This article is to propose a Deep Learning Galerkin Method (DGM) for the closed-loop geothermal system, which is a new coupled multi-physics PDEs and mainly consists of a framework of underground heat exchange pipelines to extract the geothermal heat from the geothermal reservoir. This method is a natural combination of Galerkin Method and machine learning with the solution approximated by a neural network instead of a linear combination of basis functions. We train the neural network by randomly sampling the spatiotemporal points and minimize loss function to satisfy the differential operators, initial condition, boundary and interface conditions. Moreover, the approximate ability of the neural network is proved by the convergence of the loss function and the convergence of the neural network to the exact solution in $L^2$ norm under certain conditions. Finally, some numerical examples are carried out to demonstrate the approximation ability of the neural networks intuitively.

\noindent{\bf Keywords:~~}Deep Learning Galerkin Method, deep neural network, closed-loop geothermal system, convergence, numerical experiments.

\section{Introduction}
The PDEs are used to describe a variety of physical phenomena such as fluid dynamics, quantum mechanics, and elasticity. Many subsurface energy and environmental applications are controlled by multiphase fluid flow in porous media, including hydrocarbon recovery from oil and gas reservoirs [3,8], geological carbon dioxide storage [4,6,7], geothermal energy extraction [14,18,19]. The management of such processes has been increasingly relying on high fidelity physics-based numerical simulations. As one main form of renewable green energy, geothermal energy is environmental, independent of the weather or climate condition and widely distributed, hence it can significantly contribute to the sustainable energy application [13,5,9]. The closed-loop heat exchanger, one  major technique to extract the geothermal energy, circulates working fluid down to the rock mass with high temperatures and then back to the surface through a continuous and closed-loop pipe [10,12,17,20]. There has been some studies shows that working fluid in the closed-loop geothermal system can be water, superheated steam, carbon dioxide, and so on [11,16,17].

In this article, we consider the closed-loop system for the geothermal reservoir with porous media flows [15,21,22,23]. This model combines the governing equations of the fluid flow in the pipe, the porous media flow in the geothermal reservoir, and the heat flow into a coupled system through the heat exchanging conditions and the no-fluid-communication conditions on the interface between the pipe and porous media. The proposed model is a multi-physics interface problem. Traditionally, the numerical simulators discretize the governing equations of multi-physics interface problem by finite difference/volume/element methods [26,27,28,29,30,31,32] resulting in large scale algebraic system of equations that may require significant computational resources and advanced numerical solvers and preconditioners [24,25].

During the past decade significant advances have been made in deep learning approaches [2] in both theory and application. With recent advances in CPU and GPU-based parallel processing, deep learning has achieved great success in numerous applications including image recognition, natural language processing, and automated driving to name a few. As a result, many mathematicians have introduced neural networks into the PDEs precisely because multilayer feedforward networks are proven to be universal approximators for the PDEs [33,34]. More specifically, once the network structure is determined, any order derivatives of the neural network can be obtained analytically. Coupled with the automatic differentiation technique, some studies have shown that neural networks can be applied to solve the PDEs [35-51]. Besides, many effective algorithms are proposed to solve some high-dimensional PDEs in [58-65]. Sirignano and Spiliopoulos [36] proposed the DGM to address the curse of dimensionality problem when solving high-dimensional PDEs. There are also some significant studies, Karpatne et al.[52] introduced Physics-guided Neural Networks, and successfully applied it in lake temperature modeling. Wang et al.\cite{53} developed the Theory-guided Neural Network to reliably predict subsurface single-phase flow in porous media. Yan et al.[54] described a gradient-based deep neural network constrained by the physics related to multiphase flow in porous media.

Inspired by the framework of \cite{36}, we have successfully applied the DGM for solving the second-order linear elliptic equations and stokes problems in \cite{55,56,57}. In these studies, there is few about coupled multi-physics model and multi-physics interface problem. Facing the challenge, the focus of this work is to explore the use of DNN as a universal approximator to learn the porous media flow, fluid flow and heat flow in the closed-loop system. Our algorithm is meshfree and parallel for the equations and variables, the neural network is trained on batches of randomly sampled time and space points to satisfy the differential operators, initial condition, boundary and interface conditions. Moreover, we obtain the convergence of the loss function and the neural network, in addition, the performance of the method is demonstrated by some numerical experiments.

The remaining of this paper is organized as follows. In the next section, we discuss the new coupled multi-physics model and the DGM algorithm. Section 3 is devoted to illustrate the convergence of the loss function, while we present the proof of another theorem to guarantee the convergence of neural network's solution in Section 4. Finally, some numerical results are provided in Section 5.
\section{Methodology}
\subsection{ Equations and preliminaries}
\begin{figure}[h]
  \centering
  \includegraphics[scale=0.5]{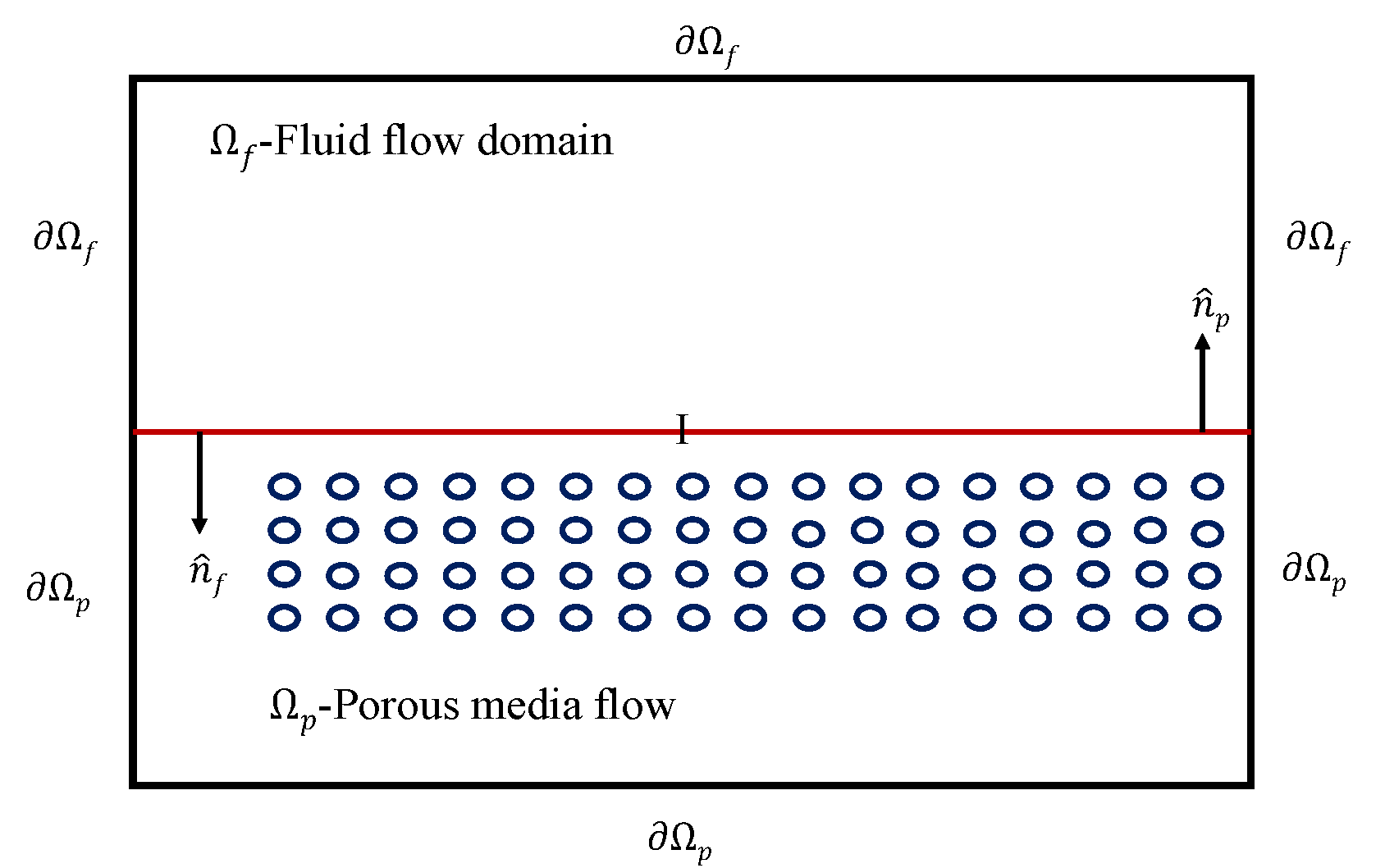}\\
  \caption{The global domain $\Omega$ separated by a common interface $\mathbb{I}$}\label{figure-domain}
\end{figure}

As showed in Figure \ref{figure-domain}, the global domain $\Omega$ consists of two subdomains $\Omega_f$ and $\Omega_p$, where $\Omega_f$, $\Omega_p\in\mathbb{R}^2$ are open, regular, simply connected, and bounded by Lipschitz continuous boundaries of $\partial\Omega_f\setminus\mathbb{I}$ and $\partial\Omega_p\setminus\mathbb{I}$. Here $\Omega_f\bigcap\Omega_p=\emptyset$, $\bar{\Omega}_f\bigcap\bar{\Omega}_p=\mathbb{I}$ and $\bar{\Omega}_f\bigcup\bar{\Omega}_p=\bar{\Omega}$. $\hat{n}_f$ and $\hat{n}_p$ are the unit normal vectors which satisfy the condition of $\hat{n}_f=-\hat{n}_p$ on the interface $\mathbb{I}$. They point outward from the free flow region $\Omega_f$ and porous media flow region $\Omega_p$ respectively. The time frame is considered in $[0,T]$.

We assume that the fluid flow with heat transfer is governed by the combination of the Navier-Stokes equations with the heat equations in $\Omega_f$:

\begin{align}
\frac{\partial u_f}{\partial t} - Pr\Delta u_f + (u_f\cdot\nabla)u_f + \nabla p_f - PrRa\xi\theta_f &= f_f \quad in~\Omega_f\times(0,T],\label{equation-2.1}\\
\nabla\cdot u_f &= 0 \quad in~\Omega_f\times(0,T],\\
u_f &= 0 \quad on~\partial\Omega_f\setminus\mathbb{I}\times(0,T],\\
u_f(0,x) &= u_f^0(x) \quad in~\Omega_f,\\
\frac{\partial\theta_f}{\partial t} - k_f\Delta\theta_f + u_f\cdot\nabla\theta_f &= \Upsilon_f \quad in~\Omega_f\times(0,T],\\
\theta_f &= 0 \quad on~\Gamma_M\times(0,T],\\
\frac{\partial\theta_f}{\partial n_f} &= 0 \quad on~\Gamma_E\times(0,T],\\
\theta_f(0,x) &= \theta_f^0(x) \quad in~\Omega_f.\label{equation-2.8}
\end{align}

Here $u_f$, $p_f$, $\theta_f$ denotes the free fluid flow region velocity vector field, pressure, and temperature respectively. The unit vector $\xi$ denote the direction of the gravitational acceleration. $f_f$ and $\Upsilon_f$ and the external force terms. $Pr$ is the Prandtl number, Ra represents the Rayleigh number, $k_f>0$ is the thermal conductivity. $\Gamma_M$, $\Gamma_E$ denotes the Dirichlet and Neumann boundary condition respectively in the pipe region boundaries where $\Gamma_M\bigcup\Gamma_E=\partial\Omega_f\setminus\mathbb{I}$.

We assume the porous media region $\Omega_p$ is homogeneous and isotropic. The porous media flow with heat transfer can be governed by the following Darcy's law coupled with heat equations:

\begin{align}
\frac{\partial u_p}{\partial t} + u_p + Da\nabla p_p - DaRa\xi\theta_p &= 0 \quad in~\Omega_p\times(0,T],\label{equation-2.9}\\
\nabla\cdot u_p &= 0 \quad in~\Omega_p\times(0,T],\\
u_p(0,x) &= u_p^0(x) \quad in~\Omega_p,\\
u_p\cdot n_p &= 0 \quad on~\partial\Omega_p\setminus\mathbb{I},\\
\frac{\partial\theta_p}{\partial t} - k_p\Delta\theta_p + u_p\cdot\nabla\theta_p &= \Upsilon_p \quad in~\Omega_p\times(0,T],\\
\theta_p &= 0 \quad on~\Gamma_N\times(0,T],\\
\frac{\partial\theta_p}{\partial n_p} &= 0 \quad on~\Gamma_Z\times(0,T],\\
\theta_p(0,x) &= \theta_p^0(x) \quad in~\Omega_p.\label{equation-2.16}
\end{align}

Here $u_p$, $p_p$ and $\theta_p$ denote the porous medium fluid flow region velocity vector field, pressure, and temperature respectively. The source term is denoted by $\Upsilon_p$, $Da$ is the Darcy number which represents the relative effect of the permeability of the medium versus the cross-sectional area, $k_p$ is the thermal conductivity. $\Gamma_N$, $\Gamma_Z$ denote the Dirichlet and Neumann boundary condition respectively in the porous media region boundaries where $\Gamma_N\bigcup\Gamma_Z=\partial\Omega_p\setminus\mathbb{I}$.

Since the closed-loop system does not have fluid communication but only heat transfer between the reservoir and the pipe, we use the following two heat exchanging interface conditions and two no-fluid-communication conditions on the interface $\mathbb{I}$.

\begin{itemize}
  \item Continuity of temperature across the interface:
  \begin{align}
  \theta_f = \theta_p~~~~~~~~\text{on}~\mathbb{I},\label{equation-2.17}
  \end{align}
  \item Continuity of heat flux across the interface:
  \begin{align}
  \hat{n}_f\cdot k_f\nabla\theta_f = - \hat{n}_p\cdot k_p\nabla\theta_p~~~~~~~~\text{on}~\mathbb{I},
  \end{align}
  \item No-communication and no-slip condition for the free flow on the interface:
  \begin{align}
  u_f\cdot\hat{n}_f=0,~~~u_f\cdot\vec{\tau}_f=0~~~~~~~~\text{on}~\mathbb{I},
  \end{align}
  where $\tau$ denote the unit tangential vector along the interface.
  \item No-porous media flow passing through the interface:
  \begin{align}
  u_p\cdot\hat{n}_p=0~~~~~~~~\text{on}~\mathbb{I}.\label{equation-2.20}
  \end{align}
\end{itemize}

In order to derive the variational formulation for the model problem \eqref{equation-2.1}-\eqref{equation-2.16}, we define the following spaces as follows:

\begin{align*}
&Y_f:=H_0^1(\Omega_f)^2:=\{v_f\in H^1(\Omega_f)^2:v_f=0~\text{on}~\partial\Omega_f\},\\
&Y_p:=H(div:\Omega_p):=\{v_p\in L^2(\Omega_p)^2,\nabla\cdot v_p\in L^2(\Omega_p):v_p\cdot\hat{n}_p=0~\text{on}~\partial\Omega_p\},\\
&T_f:=H_0^1(\Omega_f):=\{t_f\in H^1(\Omega_f):t_f=0~\text{on}~\Gamma_M\},\\
&T_p:=H_0^1(\Omega_p):=\{t_p\in H^1(\Omega_p):t_p=0~\text{on}~\Gamma_N\},\\
&Q_f:=L_0^2(\Omega_f):=\{q_f\in L^2(\omega_f):\int_{\Omega_f}q_fdx=0\},\\
&Q_p:=L_0^2(\Omega_p):=\{q_p\in L^2(\omega_p):\int_{\Omega_p}q_p dx=0\}\\
\end{align*}

Define a product space

$$W_T=T_f\times T_p,$$

and
$$W_{\mathbb{I}}=\{(t_f,t_p)\in T_f\times T_p:t_f=t_p~on~\mathbb{I}\}\subset W_T.$$

We also define the solenoidal spaces
$$V_f:=\{v_f\in Y_f:\nabla\cdot v_f=0\}~\text{and}~V_p:=\{v_p\in Y_p:\nabla\cdot v_p=0\}.$$

The trilinear form is defined as follows:
\begin{align}
c_1(u_f,v_f,t_f)_{\Omega_f}&=((u_f\cdot\nabla)v_f,t_f)+\frac{1}{2}((div u_f)v_f,t_f)\nonumber\\
&=\frac{1}{2}((u_f\cdot\nabla)v_f,t_f)-\frac{1}{2}((u_f\cdot\nabla)t_f,v_f),~~\forall u_f,v_f,t_f\in Y_f.
\end{align}

In a similar manner, we define another two trilinear forms for any  $(t_f,t_p)\in W_{\mathbb{I}}$ and $(\theta_f,\theta_p)\in W_T$:

\begin{align}
&c_2(u_f,\theta_f,t_f)_{\Omega_f}=\frac{1}{2}((u_f\cdot\nabla)\theta_f,t_f)-\frac{1}{2}((u_f\cdot\nabla)t_f,\theta_f)~~\forall u_f\in Y_f,\nonumber\\
&c_3(u_p,\theta_p,t_p)_{\Omega_p}=\frac{1}{2}((u_p\cdot\nabla)\theta_p,t_p)-\frac{1}{2}((u_p\cdot\nabla)t_p,\theta_p)~~\forall u_p\in Y_p.
\end{align}

Moreover, if $u_f\in V_f$ and $u_p\in V_p$, then for any $(u_f,v_f)\in Y_f$, $(\theta_f,\theta_p)\in W_T$, and $(t_f,t_p)\in W_{\mathbb{I}}$, we have
$$c_1(u_f,v_f,t_p)_{\Omega_f}=((u_f\cdot\nabla)v_f,t_p)_{\Omega_f},$$
$$c_2(u_f,\theta_f,t_f)_{\Omega_f}=((u_f\cdot\nabla)\theta_f,t_f)_{\Omega_f},$$
$$c_3(u_p,\theta_p,t_p)_{\Omega_p}=((u_p\cdot\nabla)\theta_p,t_p)_{\Omega_p}.$$

For simplicity of notations, let us denote
$$\mathcal{G}_f^1[u_f,p_f,\theta_f] = \frac{\partial u_f}{\partial t} - Pr\Delta u_f + (u_f\cdot\nabla)u_f + \nabla p_f - PrRa\xi\theta_f  - f_f,$$
$$\mathcal{G}_f^2 [u_f,\theta_f]= \frac{\partial\theta_f}{\partial t} - k_f\Delta\theta_f + u_f\cdot\nabla\theta_f - \Upsilon_f,$$
$$\mathcal{G}_p^1[u_p,p_p,\theta_p] = \frac{\partial u_p}{\partial t} + u_p + Da\nabla\phi_p - DaRa\xi\theta_p,$$
$$\mathcal{G}_p^2[u_p,\theta_p] = \frac{\partial\theta_p}{\partial t} - k_p\Delta\theta_p + u_p\cdot\nabla\theta_p - \Upsilon_p,$$
where the $(u_f,p_f,\theta_f)$ is the analytic solution of the system \eqref{equation-2.1}-\eqref{equation-2.8} and the $(u_p,p_p,\theta_p)$ is the analytic solution of the system \eqref{equation-2.9}-\eqref{equation-2.16}.

\subsection{The DGM algorithm}

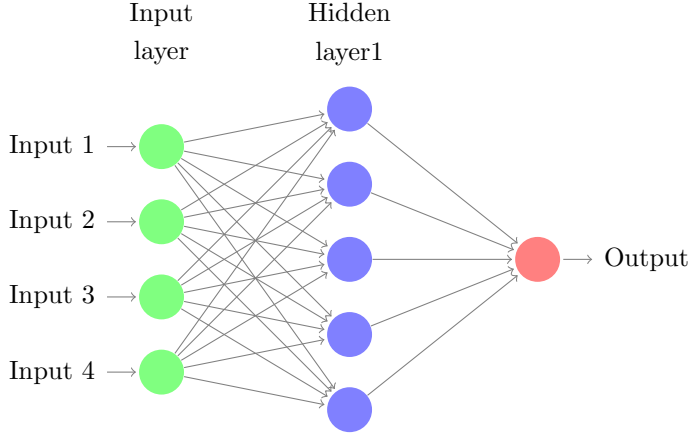
\begin{figure}[H]
\centering
\def\layersep{2.5cm}

  \begin{tikzpicture}[shorten >=1pt,->,draw=black!50, node distance=\layersep]
     \tikzstyle{every pin edge}=[<-,shorten <=1pt]
     \tikzstyle{neuron}=[circle,fill=black!25,minimum size=17pt,inner sep=0pt]
     \tikzstyle{input neuron}=[neuron, fill=green!50];
     \tikzstyle{output neuron}=[neuron, fill=red!50];

    \tikzstyle{hidden neuron1}=[neuron, fill=blue!50];

     \tikzstyle{annot} = [text width=4em, text centered]

     \foreach \name / \y in {1,...,4}
        \node[input neuron, pin=left:Input \y] (I-\name) at (0,-\y) {};

     \foreach \name / \y in {1,...,5}
         \path[yshift=0.5cm]
             node[hidden neuron1] (H1-\name) at (\layersep,-\y cm) {};

    \node[output neuron,pin={[pin edge={->}]right:Output}, right of=H1-3] (O) {};

     \foreach \source in {1,...,4}
         \foreach \dest in {1,...,5}
             \path (I-\source) edge (H1-\dest);
     %
     \foreach \source in {1,...,5}
         \path (H1-\source) edge (O);

     \node[annot,above of=H1-1, node distance=1cm] (hl) {Hidden layer1};
     \node[annot,left of=hl] {Input layer};
 \end{tikzpicture}
\caption{The neural networks with one hidden layer and one output unit.}\label{figure-NN}
\end{figure}

We apply the multilayer feed forward networks, which consists of an input layer, one or more hidden layers and one output layer, to approximate solution of the equations \eqref{equation-2.1}-\eqref{equation-2.16}. The approximation capabilities of neural network architectures have recently been investigated by many authors. Especially, in \cite{60}, they have ascertained that the standard multilayer feedforward networks with activation function $\psi$ can approximate any continuous function well on arbitrary compact subsets $X$ of $\Omega$, whenever activation function $\psi$ is continuous, bounded and nonconstant. For convenience, we will specifically contrive our results only for the instance where there is only one hidden layer, each hidden unit has the same activation function $\psi_i$ and one output unit. The simplified diagram of this neural network is shown in Figure \ref{figure-NN}.

The set of all the functions are described by such a network with $n$ hidden units for velocity, pressure and temperature of $\Omega_f$ and $\Omega_p$, respectively.

$$[\mathfrak{C}_{u_f}^n(\varphi_f)]^d=\Big\{\Phi_f(x,t):R^{d+1} \rightarrow R^d|\Phi_f^l(x,t)=\sum_{i=1}^n\beta_i^1\varphi_f(\sigma^1_{1,i}t+\sum_{j=1}^d\sigma_{j,i}^1x_j+c_{1,i})\Big\},$$
$$[\mathfrak{C}_{u_p}^n(\varphi_p)]^d=\Big\{\Phi_p(x,t):R^{d+1} \rightarrow R^d|\Phi_p^l(x,t)=\sum_{i=1}^n\beta_i^2\varphi_p(\sigma^2_{1,i}t+\sum_{j=1}^d\sigma_{j,i}^2x_j+c_{2,i})\Big\},$$
$$\mathfrak{C}_{p_f}^n(\psi_f)=\Big\{\Psi_f(x,t):R^{d+1} \rightarrow R|\Psi_f(x,t)=\sum_{i=1}^n\beta_i^3\psi_f(\sigma^3_{1,i}t+\sum_{j=1}^d\sigma_{j,i}^3x_j+c_{3,i})\Big\},$$
$$\mathfrak{C}_{p_p}^n(\psi_p)=\Big\{\Psi_p(x,t):R^{d+1} \rightarrow R|\Psi_p(x,t)=\sum_{i=1}^n\beta_i^4\psi_p(\sigma^4_{1,i}t+\sum_{j=1}^d\sigma_{j,i}^4x_j+c_{4,i})\Big\},$$
$$\mathfrak{C}_{\theta_f}^n(\omega_f)=\Big\{\Theta_f(x,t):R^{d+1} \rightarrow R|\Theta_f(x,t)=\sum_{i=1}^n\beta_i^5\omega_f(\sigma^5_{1,i}t+\sum_{j=1}^d\sigma_{j,i}^5x_j+c_{5,i})\Big\},$$
and
$$\mathfrak{C}_{\theta_p}^n(\omega_p)=\Big\{\Theta_p(x,t):R^{d+1} \rightarrow R|\Theta_p(x,t)=\sum_{i=1}^n\beta_i^6\omega_p(\sigma^6_{1,i}t+\sum_{j=1}^d\sigma_{j,i}^6x_j+c_{6,i}) \Big\}.$$
Especially,
$$\mathfrak{C}_{\theta_\mathbb{I}}^n(\omega_{\mathbb{I}})=\Big\{\Big(\Theta_f(x,t),\Theta_p(x,t)\Big)\in\mathfrak{C}_{\theta}^n(\omega_f)\times\mathfrak{C}_{\theta}^n(\omega_p):\Theta_f(x,t)|_{\mathbb{I}}= \Theta_p(x,t)|_{\mathbb{I}}\Big\},$$
is a product space.

\begin{itemize}
  \item $\mathfrak{C}^{n}_{*}$: the class of neural networks consisting of a single hidden layer with $n$ hidden units and $d$ random points $x_j,~j=1,2,\cdots, d$ and one output unit, for convenience, we still use this notation to represent multilayer neural networks;
  \item $\varphi_i,\psi_i,\omega_i,(i=f,p)$ are shared activate functions of the hidden units, bounded and non-constant.
  \item $\beta_{i}^{k},\sigma_{j}^{k},\sigma_{ji}^{k}$ are weights, $c_{k,i}$ are the thresholds of the network.
  \item Denote $\Lambda_k=(\beta_1^k,\cdots,\beta_n^k,\sigma_{1,1}^k,\sigma_{1,2}^K,\cdots,\sigma_{1,n}^k,\sigma_{1,1}^k,\cdots,\sigma_{d,n}^k,c_{k,1},\cdots,c_{k,n})\in R^{2n+(1+d)n},(k=1,\cdots,6)$ are the parameters of the neural network.
\end{itemize}

Assuming that $(U_f,P_f,\Theta_f)$ and $(U_p,P_p,\Theta_p)$ are neural network solutions to the system \eqref{equation-2.1}-\eqref{equation-2.8} and \eqref{equation-2.9}-\eqref{equation-2.16}, respectively. Then the loss function of the two systems are defined by
\begin{align}
J_f=&\Big\|\mathcal{G}_f^1[U_f,P_f,\Theta_f](x,t;\Lambda_1,\Lambda_3,\Lambda_5)\Big\|^2_{0,\Omega_{f},\nu_1}+\Big\|\mathcal{G}_f^2 [U_f,\Theta_f](x,t;\Lambda_1,\Lambda_5)\Big\|^2_{0,\Omega_{f},\nu_1}\nonumber\\
&+\Big\|\nabla\cdot U_f(x,t;\Lambda_1)\Big\|^2_{0,\Omega_{f},\nu_1}+\Big\|U_f(x,t;\Lambda_1)\Big\|^2_{0,\partial\Omega_f\setminus\mathbb{I},\nu_2}+\Big\|U_f(x,0;\Lambda_1)-u_f^0(x)\Big\|^2_{0,\Omega_f,\nu_3}\nonumber\\
&+\Big\|\Theta_f(x,t;\Lambda_5)\Big\|^2_{0,\Gamma_M,\nu_2}+\Big\|\Theta_f(x,0;\Lambda_5)-\theta_f^0(x)\Big\|^2_{0,\Omega_f,\nu_3} +\Big\|\Theta_f\cdot
\hat{n}_p(x,t;\Lambda_5)\Big\|^2_{0,\Gamma_E,\nu_2},\\
J_p=&\Big\|\mathcal{G}_p^1[U_p,P_p,\Theta_p](x,t;\Lambda_2)\Big\|^2_{0,\Omega_{p},\nu_4}+\Big\|\mathcal{G}_p^2[U_p,\Theta_p](x,t;\Lambda_2,\Lambda_6)\Big\|^2_{0,\Omega_{p},\nu_4}\nonumber\\
&+\Big\|\nabla\cdot U_p(x,t;\Lambda_2)\Big\|^2_{0,\Omega_{p},\nu_4}+\Big\|U_p(x,0;\Lambda_2)-u_p^0(x)\Big\|^2_{0,\Omega_p,\nu_5}\nonumber\\
&+\Big\|\Theta_p(x,t;\Lambda_6)\Big\|^2_{0,\Gamma_N,\nu_6}+\Big\|\Theta_p(x,0;\Lambda_6)-\theta_p^0(x)\Big\|^2_{0,\Omega_p,\nu_5}\nonumber\\
&+\Big\|U_p\cdot \hat{n}_p(x,t;\Lambda_6)\Big\|^2_{0,\partial\Omega_p\setminus\mathbb{I},\nu_6}+\Big\|\Theta_p\cdot \hat{n}_p(x,t;\Lambda_6)\Big\|^2_{0,\Gamma_Z,\nu_6},\\
\end{align}
and
\begin{align}
J_{\mathbb{I}}=&\Big\|(\Theta_f-\Theta_p)(x,t;\Lambda_5,\Lambda_6)\Big\|^2_{0,\mathbb{I},\nu_7}+\Big\|(U_f\cdot \hat{n}_f)(x,t;\Lambda_1)\Big\|^2_{0,\mathbb{I},\nu_7}+\Big\|(U_p\cdot \hat{n}_p)(x,t;\Lambda_2)\Big\|^2_{0,\mathbb{I},\nu_7}\nonumber\\
&+\Big\|(\hat{n}_f\cdot k_f\nabla\Theta_f+\hat{n}_p\cdot k_p\nabla\Theta_p)(x,t;\Lambda_5,\Lambda_6)\Big\|^2_{0,\mathbb{I},\nu_7},\\
J =& J_f+J_p+J_{\mathbb{I}}.\label{equation-2.26}
\end{align}

There are some notes for the loss function:
\begin{itemize}
  \item $\|f(y)\|^2_{0,\mathcal{Y},\nu}=\int_{\mathcal{Y}}|f(y)|^2\nu(y)dy$, where $\nu$ is a positive probability density on $y\in \mathcal{Y}$;
  \item $\Lambda=\{\Lambda_1,\Lambda_2,\Lambda_3,\Lambda_4,\Lambda_5,\Lambda_6\}$ are the parameters of the neural network, $\nu_{i}(i=1,2,3)$  are the positive probability densities of $(\Omega_{f},\Omega_{p})$, $\partial\Omega\setminus\mathbb{I}$ and $\mathbb{I}$, respectively;
  \item $J_i(i=f,p,\mathbb{I})$ is defined to measure the difference between the output of the network with the known input. The whole deep learning process is to find appropriate parameters $\Lambda$ such that $J$ is as close to 0 as possible.
\end{itemize}

The DGM is provided as follows:

{\bf Step 1.}~~Generate sample points ${x_{n}, t_{n}}\in\Omega\times[0,T]$ and ${z_{n},\tau_{n}}\in\partial\Omega\times[0,T]$ based on respective probability densities $\nu_{i}(i=1,2,3)$ randomly.

{\bf Step 2.}~~At the random points $s_n=\{(x_n,t_n),(z_n,\tau_n)\}$, we compute the square error $G(\Lambda_n;s_n)$ by
$$G(\Lambda_n;s_n)=G_f(\Lambda_n;s_n)+G_p(\Lambda_n;s_n)+G_{\mathbb{I}}(\Lambda_n;s_n).$$

{\bf Step 3.}~~Apply the stochastic gradient descent step at $s_{n}$ by
$$\Lambda_{n+1}=\Lambda_n-\eta_n\nabla_{\Lambda}G(\Lambda_n;s_n).$$
Here, the learning rate $\eta_n$ decreases as $n$ increases.

{\bf Step 4.}~~Repeat step 3 until $G(\Lambda_n;s_n)$ meets the convergence condition $$\lim_{n\rightarrow\infty}\nabla_{\Lambda}G(\Lambda_n,s_n)=0,$$
where $\mathbb{E}\Big[\nabla_{\Lambda}G(\Lambda_n,s_n)|_{\Lambda_n}\Big]=\nabla_{\Lambda}J(\Lambda_n,s_n)$, i.e. $\nabla_{\Lambda}G(\Lambda_n,s_n)|_{\Lambda_n}$ is the unbiased estimation of $\nabla_{\Lambda}J(\Lambda_n,s_n)$.
\begin{algorithm}
        \caption{ Deep Learning Galerkin Method.}
        \begin{algorithmic}[1] 
            \Require training data, $s_n=\{(x_n,t_n),(z_n,\tau_n)\}$, max number of iterations M;
            \Ensure Neural network predicted solutions, $U_f(x,t)$, $U_p(x,t)$, $\Theta_f(x,t)$, $\Theta_p(x,t)$, $P_f(x,t)$, $P_p(x,t)$;
            \State Specify the training set $s_n=\{(x_n,t_n),(z_n,\tau_n)\}$ for the initial/boundary conditions, domain and equations;
            \State Construct a neural network $\{U_f(x,t;\Lambda), U_p(x,t;\Lambda), \Theta_f(x,t;\Lambda), \Theta_p(x,t;\Lambda), P_f(x,t;\Lambda), P_p(x,t;\Lambda)\}$;
            \State Specify a loss function by summing the initial/boundary conditions, the residual of the PDE;
            \State Train the neural network to find the best parameters $\Lambda$ by minimizing the loss function $J$ in \eqref{equation-2.26};
            \State Get the predicted solution on the domain;\\
            \Return Return $U_f(x,t)$, $U_p(x,t)$, $\Theta_f(x,t)$, $\Theta_p(x,t)$, $P_f(x,t)$, $P_p(x,t)$.

        \end{algorithmic}
\end{algorithm}

\section{Convergence rates}

\subsection{Convergence of the loss function}

{$\textbf{\emph{(H1)}}$.}~~Let the terms $u_f(x)$ and $\Delta u_f(x)$ be locally $Lipschitz$ in $(u_f,\nabla u_f)$ with $Lipschitz$ constant which has at most polynomial growth on $u$, uniformly with respect to $x$. This means that
$$|\Delta U_f-\Delta u_f|\leq\Big(|\nabla U_f|^{q_1/2}+|\nabla u_f|^{q_2/2}\Big)|\nabla U_f-\nabla u_f|,~~0\leq q_i\leq\infty,~~i=1,2.$$

\begin{align*}
|\Delta U_f-\Delta u_f|^2&\leq\Big(|\nabla U_f|^{q_1/2}+|\nabla u_f|^{q_2/2}\Big)^2|\nabla U_f-\nabla u_f|^2\\
&\leq\Big(|\nabla U_f|^{q_1}+|\nabla u_f|^{q_2}+2|\nabla U_f|^{q_1/2}\cdot|\nabla u_f|^{q_2/2}\Big)|\nabla U_f-\nabla u_f|^2\\
&\leq2\Big(|\nabla U_f|^{q_1}+|\nabla u_f|^{q_2}\Big)|\nabla U_f-\nabla u_f|^2.
\end{align*}
For the $(U_p,u_p)$, $(\Theta_f,\theta_f)$, $(\Theta_p,\theta_p)$, $(P_f,p_f)$, $(P_p,p_p)$, the above properties are utilized.

{$\textbf{\emph{Lemma~1}}$.}~~If $u_f\in Y_f$, then there exists a positive constant $C>0$ such that
$$\|u_f\|_{L^4(\Omega_f)}\leq C\|u_f\|_{L_2(\Omega_f)}^{\frac{1}{2}}\|\nabla u_f\|_{L_2(\Omega_f)}^{\frac{1}{2}}$$

{$\textbf{\emph{Lemma~2}}$.}~~The weak formulation of the closed-loop geothermal system \eqref{equation-2.1}-\eqref{equation-2.16} has a unique weak solution $(u_f,u_p;p_f,p_p;\theta_f,\theta_p)\in(Y_f\times Y_p\times Q_f\times Q_p\times W_T)$ satisfying
$$u_f\in L^{\infty}(0,T;L^2(\Omega_f))\cap L^2(0,T;Y_f),~u_p\in L^{\infty}(0,T;L^2(\Omega_p))\cap L^2(0,T;L^2(\Omega_p)),~p_f\in L^2(0,T;Q_f),$$
$$p_p\in L^2(0,T;Q_f),~\theta_f\in L^{\infty}(0,T;L^2(\Omega_f))\cap L^2(0,T;T_f),~\theta_p\in L^{\infty}(0,T;L^2(\Omega_p))\cap L^2(0,T;T_p).$$

{$\textbf{\emph{Lemma~3}}$.}~~By Theorem 3 of \cite{33} we know that there exist functions
\begin{itemize}
  \item $U_f\in[\mathfrak{C}_{u_f}^n(\varphi_f)]^d$, $U_p\in[\mathfrak{C}_{u_p}^n(\varphi_p)]^d$ are uniformly 2-dense on compact sets of $\mathcal{C}^{\infty}(0,T;L^2(\Omega_f))$ and $\mathcal{C}^{\infty}((0,T;L^2(\Omega_p))$, respectively. This means that for $u_f\in\mathcal{C}^{\infty}(0,T;L^2(\Omega_f))$, $u_p\in\mathcal{C}^{\infty}(0,T;L^2(\Omega_p)$ and $\epsilon_1>0$, $\epsilon_2>0$, there exist
      $$\sup_{(t,x)\in\Omega_f}|\partial_t U_f(t,x;\Lambda_1)-\partial_t u_f(t,x)|+\max_{|\alpha|\leq2}\sup_{x\in{\Omega_f}}|\partial_{x}^{(\alpha)}U_f(x,t;\Lambda_1)-\partial_{x}^{(\alpha)}u_f(x,t)|<\epsilon_1,$$
      $$\sup_{(t,x)\in\Omega_p}|\partial_t U_p(t,x;\Lambda_2)-\partial_t
      u_p(t,x)|+\max_{|\alpha|\leq2}\sup_{x\in{\Omega_p}}|\partial_{x}^{(\alpha)}U_p(x,t;\Lambda_2)-\partial_{x}^{(\alpha)}u_p(x,t)|<\epsilon_2;$$
  \item $\Theta_f$ and $\Theta_p\in\mathfrak{C}_{\theta_\mathbb{I}}^n(\omega_{\mathbb{I}})$ are uniformly 2-dense on compact sets of $\mathcal{C}^{\infty}(0,T;L^2(\Omega_f)\times L^2(\Omega_p))$, the meaning is that for $\theta_f,\theta_p\in\mathcal{C}^{\infty}(0,T;L^2(\Omega_f)\times L^2(\Omega_p))$ and $\epsilon_3>0$, $\epsilon_4>0$, there exist
      $$\sup_{(t,x)\in\Omega_f}|\partial_t \Theta_f(t,x;\Lambda_5)-\partial_t \theta_f(t,x)|+\max_{|\alpha|\leq2}\sup_{x\in{\Omega_f}}|\partial_{x}^{(\alpha)}\Theta_f(x,t;\Lambda_5)-\partial_{x}^{(\alpha)}\theta_f(x,t)|<\epsilon_3,$$
      $$\sup_{(t,x)\in\Omega_p}|\partial_t \Theta_p(t,x;\Lambda_6)-\partial_t \theta_p(t,x)|+\max_{|\alpha|\leq2}\sup_{x\in{\Omega_p}}|\partial_{x}^{(\alpha)}\Theta_p(x,t;\Lambda_6)-\partial_{x}^{(\alpha)}\theta_p(x,t)|<\epsilon_4;$$
  \item $P_f\in\mathfrak{C}_{p_f}^n(\psi_f)$, $P_p\in\mathfrak{C}_{p_p}^n(\psi_p)$ are uniformly 2-dense on compact sets of $\mathcal{C}^{2}(0,T;L^2(\Omega_f))$ and $\mathcal{C}^{2}(0,T;L^2(\Omega_p))$, respectively. This means that for $p_f\in\mathcal{C}^{2}(0,T;L^2(\Omega_f))$, $p_p\in\mathcal{C}^{2}(0,T;L^2(\Omega_p))$ and $\epsilon_5>0$, $\epsilon_6>0$, there exist
        $$\sup_{(t,x)\in\Omega_f}|\partial_t P_f(t,x;\Lambda_3)-\partial_t p_f(t,x)|<\epsilon_5,$$
        $$\sup_{(t,x)\in\Omega_p}|\partial_t P_p(t,x;\Lambda_4)-\partial_t p_p(t,x)|<\epsilon_6.$$
\end{itemize}

\noindent{\bf Theorem 1.~~}Under the above assumptions, there exists a positive constant $K>0$ and a neural network solution $\overline{U}\in \Big([\mathfrak{C}_{u_f}^n(\varphi_f)]^d\times[\mathfrak{C}_{u_p}^n(\varphi_p)]^d\times\mathfrak{C}_{p_f}^n(\psi_f)\times\mathfrak{C}_{p_p}^n(\psi_p)\times\mathfrak{C}_{\theta_\mathbb{I}}^n(\omega_{\mathbb{I}})\Big)$ such that, for $\forall\epsilon>0$,
$$J(\overline{U})\leq K\epsilon,$$
where $\overline{U}= (U_f,U_p,P_f,P_p,\Theta_f,\Theta_p)$, $K$ depends on $f_f$, $\Upsilon_f$, $\Upsilon_p$ and $\nu_{i}(i=1,2,3)$.

$Proof.~~$Let $u_f$, $p_f$, $\theta_f$ be the solution to \eqref{equation-2.1}-\eqref{equation-2.8}, then $\mathcal{G}_f^1[u_f,p_f,\theta_f](x)=0$ naturally. Thus, we have

\begin{align}
J_f=&\Big\|\mathcal{G}_f^1[U_f,P_f,\Theta_f]-\mathcal{G}_f^1[u_f,p_f,\theta_f]\Big\|^2_{0,\Omega_{f},\nu_1}+\Big\|\mathcal{G}_f^2 [U_f,\Theta_f]-\mathcal{G}_f^2 [u_f,\theta_f]\Big\|^2_{0,\Omega_{f},\nu_1}\nonumber\\
&+\Big\|\nabla\cdot (U_f-u_f)\Big\|^2_{0,\Omega_{f},\nu_1}+\Big\|U_f\Big\|^2_{0,\partial\Omega_f\setminus\mathbb{I},\nu_2}+\Big\|U_f(0,x)-u_f^0(x)\Big\|^2_{0,\Omega_f,\nu_1}\nonumber\\
&+\Big\|\Theta_f\Big\|^2_{0,\Gamma_M,\nu_2}+\Big\|\Theta_f(0,x)-\Theta_f^0(x)\Big\|^2_{0,\Omega_f,\nu_1}+\Big\|\Theta_f\cdot n_p\Big\|^2_{0,\Gamma_E,\nu_2}.
\end{align}
In the following we will estimate the terms on the right hand side of $(3.1)$ one by one as follows. First, we have
\begin{align}
&\Big\|\mathcal{G}_f^1[U_f,P_f,\Theta_f]-\mathcal{G}_f^1[u_f,p_f,\theta_f]\Big\|^2_{0,\Omega_{f},\nu_1}\nonumber\\
\leq&\int_{\Omega_{f}}\Big|\frac{\partial U_f}{\partial t}-\frac{\partial u_f}{\partial t}\Big|^2d\nu_1(t,x)+Pr\int_{\Omega_{T_f}}\Big|\Delta U_f-\Delta u_f\Big|^2d\nu_1(t,x)\nonumber\\
&+\int_{\Omega_{f}}\Big|\nabla P_f-\nabla p_f\Big|^2d\nu_1(t,x)+PrRa\xi\int_{\Omega_{f}}\Big|\Theta_f-\theta_f\Big|^2d\nu_1(t,x)\nonumber\\
&+\int_{\Omega_{f}}\Big|(U_f\cdot\nabla)U_f-(u_f\cdot\nabla)u_f\Big|^2d\nu_1(t,x)\nonumber\\
\leq&\epsilon_1^2+2Pr\epsilon_1^2\Big(\epsilon_1^{q_1}+\sup|\nabla u_f|^{q_{1}\vee q_{2}}\Big)\epsilon_1^{2}+\epsilon_5^2+PrRa\xi\epsilon_3^2,
\end{align}
by applying the hypothesis $\emph{(H1)}$ and $\emph{Lemma~3}$.\\
Similarly, we can obtain
\begin{align}
&\Big\|\mathcal{G}_f^2 [U_f,\Theta_f]-\mathcal{G}_f^2 [u_f,\theta_f]\Big\|^2_{0,\Omega_{f},\nu_1}\nonumber\\
\leq&\int_{\Omega_f}\Big|\frac{\partial\Theta_f}{\partial t}-\frac{\partial\theta_f}{\partial t}\Big|^2d\nu_1(t,x)+k_f\int_{\Omega_f}\Big|\Delta\Theta_f-\Delta\theta_f\Big|^2d\nu_1(t,x)\nonumber\\
&+\int_{\Omega_f}\Big|U_f\cdot\nabla\Theta_f-u_f\cdot\nabla\theta_f\Big|^2d\nu_1(t,x)\nonumber\\
\leq&\epsilon_3^2+2k_f\Big(\epsilon_3^{q_1}+\sup|\nabla \theta_f|^{q_{1}\vee q_{2}}\Big)\epsilon_3^{2}.
\end{align}
Add (3.3) and (3.2) to (3.1), we can get that
\begin{align}
J_f\leq&\Big\|\mathcal{G}_f^1[U_f,P_f,\Theta_f]-\mathcal{G}_f^1[u_f,p_f,\theta_f]\Big\|^2_{0,\Omega_{f},\nu_1}+\Big\|\mathcal{G}_f^2 [U_f,\Theta_f]-\mathcal{G}_f^2 [u_f,\theta_f]\Big\|^2_{0,\Omega_{f},\nu_1}\nonumber\\
&+\int_{\Omega_{f}}\Big|\nabla\cdot(U_f-u_f)\Big|^2d\nu_1(t,x)+\int_{\partial\Omega_f\setminus\mathbb{I}}\Big|U_f\Big|^2d\nu_2(t,x)+\int_{\Omega_f}\Big|U_f(0,x)-u_f^0(x)\Big|^2d\nu_1(t,x)\nonumber\\
&+\int_{\Gamma_M}\Big|\Theta_f\Big|^2d\nu_2(t,x)+\int_{\Omega_f}\Big|\Theta_f(0,x)-\Theta_f^0(x)\Big|^2d\nu_1(t,x)+\int_{\Gamma_E}\Big|\Theta_f\Big|^2d\nu_2(t,x)\nonumber\\
\leq&4\epsilon_1^2+2Pr\epsilon_1^2\Big(\epsilon_1^{q_1}+\sup|\nabla u_f|^{q_{1}\vee q_{2}}\Big)\epsilon_1^2+\epsilon_5^2+PrRa\xi\epsilon_3^2+4\epsilon_3^2+2k_f\Big(\epsilon_1^{q_1}+\sup|\nabla \theta_f|^{q_{1}\vee q_{2}}\Big)\epsilon_3^{2}\nonumber\\
\leq& K_1\epsilon^2.
\end{align}

It is similar to $J_f$,
\begin{align}
J_p=&\Big\|\mathcal{G}_p^1[U_p,P_p,\Theta_p]-\mathcal{G}_p^1[u_p,p_p,\theta_p]\Big\|^2_{0,\Omega_{p},\nu_1}+\Big\|\mathcal{G}_p^2[U_p,\Theta_p]-\mathcal{G}_p^2[u_p,\theta_p]\Big\|^2_{0,\Omega_{p},\nu_1}\nonumber\\
&+\Big\|\nabla\cdot U_p\Big\|^2_{0,\Omega_{p},\nu_1}+\Big\|U_p(0,x)-u_p^0(x)\Big\|^2_{0,\Omega_p,\nu_1}+\Big\|U_p\cdot \hat{n}_p\Big\|^2_{0,\partial\Omega_P\setminus\mathbb{I},\nu_2}\nonumber\\
&+\Big\|\Theta_p\Big\|^2_{0,\Gamma_N,\nu_2}+\Big\|\Theta_p(0,x)-\Theta_p^0(x)\Big\|^2_{0,\Omega_p,\nu_1}+\Big\|\Theta_p\cdot \hat{n}_p\Big\|^2_{0,\Gamma_Z,\nu_2}.
\end{align}
Then we estimate the right hand side of $(3.5)$ using the same method,
\begin{align}
&\Big\|\mathcal{G}_p^1[U_p,P_p,\Theta_p]-\mathcal{G}_p^1[u_p,p_p,\theta_p]\Big\|^2_{0,\Omega_{T_p},\nu_1}\nonumber\\
\leq&\int_{\Omega_{T_p}}\Big|\frac{\partial U_p}{\partial t}-\frac{\partial u_p}{\partial t}\Big|^2d\nu_1(t,x)+Da\int_{\Omega_{T_p}}\Big|\nabla P_p-\nabla p_p\Big|^2d\nu_1(t,x)\nonumber\\
&+\int_{\Omega_{T_p}}\Big|U_p-u_p\Big|^2d\nu_1(t,x)+DaRa\xi\int_{\Omega_{T_p}}\Big|\Theta_p-\theta_p\Big|^2d\nu_1(t,x)\nonumber\\
\leq&2\epsilon_2^2+Da\epsilon_6^2+DaRa\xi\epsilon_4^2,
\end{align}
and
\begin{align}
&\Big\|\mathcal{G}_p^2[U_p,\Theta_p]-\mathcal{G}_p^2[u_p,\theta_p]\Big\|^2_{0,\Omega_{T_p},\nu_1}\nonumber\\
\leq&\int_{\Omega_p}\Big|\frac{\partial\Theta_p}{\partial t}-\frac{\partial\theta_p}{\partial t}\Big|^2d\nu_1(t,x)+k_p\int_{\Omega_p}\Big|\Delta\Theta_p-\Delta\theta_p\Big|^2d\nu_1(t,x)\nonumber\\
&+\int_{\Omega_p}\Big|U_p\cdot\nabla\Theta_p-u_p\cdot\nabla\theta_p\Big|^2d\nu_1(t,x)\nonumber\\
\leq& \epsilon_4^2+\Big(\epsilon_4^{q_1}+\sup|\nabla \theta_p|^{q_{1}\vee q_{2}}\Big)\epsilon_4^{2}.
\end{align}
Plugging the above estimates into (3.5), we obtain that
\begin{align}
J_p=&\Big\|\mathcal{G}_p^1[U_p,P_p,\Theta_p]-\mathcal{G}_p^1[u_p,p_p,\theta_p]\Big\|^2_{0,\Omega_{T_p},\nu_1}+\Big\|\mathcal{G}_p^2[U_p,\Theta_p]-\mathcal{G}_p^2[u_p,\theta_p]\Big\|^2_{0,\Omega_{T_p},\nu_1}\nonumber\\
&+\int_{\Omega_{T_p}}\Big|\nabla\cdot (U_p-u_p)\Big|^2d\nu_1(t,x)+\int_{\Omega_p}\Big|U_p(0,x)-u_p^0(x)\Big|^2d\nu_1(t,x)+\int_{\partial\Omega_P\setminus\mathbb{I}}\Big|U_p\Big|^2d\nu_2(t,x)\nonumber\\
&+\int_{\Gamma_N}\Big|\Theta_p\Big|^2d\nu_2(t,x)+\int_{\Omega_p}\Big|\Theta_p(0,x)-\Theta_p^0(x)\Big|^2d\nu_1(t,x)+\int_{\Gamma_Z}\Big|\Theta_p\Big|^2d\nu_2(t,x)\nonumber\\
\leq&5\epsilon_2^2+Da\epsilon_6^2+DaRa\xi\epsilon_4^2+4\epsilon_4^2+\Big(\epsilon_4^{q_1}+\sup|\nabla \theta_p|^{q_{1}\vee q_{2}}\Big)\epsilon_4^{2}\nonumber\\
\leq& K_2\epsilon.
\end{align}

Finally, we have
\begin{align}
J_{\mathbb{I}}=&\Big\|\Theta_f-\Theta_p\Big\|^2_{0,\mathbb{I},\nu_3}+\Big\|U_f\cdot \hat{n}_f\Big\|^2_{0,\mathbb{I},\nu_3}+\Big\|U_p\cdot \hat{n}_p\Big\|^2_{0,\mathbb{I},\nu_3}+\Big\|\hat{n}_f\cdot k_f\nabla\Theta_f+\hat{n}_p\cdot k_p\nabla\Theta_p\Big\|_{0,\mathbb{I},\nu_3}\nonumber\\
\leq&\Big\|\Theta_f\Big\|^2_{0,\mathbb{I},\nu_3}+ \Big\|\Theta_p\Big\|^2_{0,\mathbb{I},\nu_3}+\Big\|U_f\Big\|^2_{0,\mathbb{I},\nu_3}+\Big\|U_p\Big\|^2_{0,\mathbb{I},\nu_3}+k_f\Big\|\nabla\Theta_f\Big\|^2_{0,\mathbb{I},\nu_3}+k_p\Big\|\nabla\Theta_p\Big\|^2_{0,\mathbb{I},\nu_3}\nonumber\\
\leq&\Big\|U_f\Big\|^2_{0,\mathbb{I},\nu_3}+\Big\|U_p\Big\|^2_{0,\mathbb{I},\nu_3}+k_f\Big\|\nabla\Theta_f\Big\|^2_{0,\mathbb{I},\nu_3}+k_p\Big\|\nabla\Theta_p\Big\|^2_{0,\mathbb{I},\nu_3}\nonumber\\
=&\Big\|U_f-u_f\Big\|^2_{0,\mathbb{I},\nu_3}+\Big\|U_p-u_p\Big\|^2_{0,\mathbb{I},\nu_3}+\Big\|u_f\Big\|^2_{0,\mathbb{I},\nu_3}+\Big\|u_p\Big\|^2_{0,\mathbb{I},\nu_3}\nonumber\\
&+k_f\Big\|\nabla(\Theta_f-\theta_f)\Big\|^2_{0,\mathbb{I},\nu_3}+k_p\Big\|\nabla(\Theta_p-\theta_p)\Big\|^2_{0,\mathbb{I},\nu_3}+k_f\Big\|\nabla\theta_f\Big\|^2_{0,\mathbb{I},\nu_3}+k_p\Big\|\nabla\theta_p\Big\|^2_{0,\mathbb{I},\nu_3}\nonumber\\
\leq&\epsilon_1^2+\epsilon_2^2+\sup(|u_f|^2+|u_p|^2)+k_f\epsilon_3^2+k_p\epsilon_4^2+k_f\sup|\nabla\theta_f|^2+k_p\sup|\nabla\theta_p|^2\nonumber\\
\leq& K_3\epsilon^2.
\end{align}
There is
$$J_f+J_p+J_\mathbb{I}\leq K\epsilon,\quad\text{i.e.},\quad J\leq K\epsilon$$
by (3.4), (3.8) and (3.9).

\subsection{Convergence of the neural network to the unique solution}
This subsection contains convergence results of neural networks to the unique solution of the closed-loop geothermal system as $n\rightarrow\infty$. During the whole parallel computing process, each variable can be calculated independently. For convenience of explanation, we take the same network
structure for each variable. Inspired by the thought of the Galerkin Method, each neural network $(U_f^n,U_p^n,P_f^n,P_p^n,\Theta_f^n,\Theta_p^n)$ satisfies the system
\begin{align}
\mathcal{G}_f^1[U_f^n,P_f^n,\Theta_f^n](x,t)&=0\quad in~\Omega_f\times(0,T],\\
\nabla\cdot U_f^n &= 0 \quad in~\Omega_f\times(0,T],\nonumber\\
U_f^n(x,t) &= 0 \quad on~\partial\Omega_f\setminus\mathbb{I}\times(0,T],\nonumber\\
U_f^n(0,x) &= U_f^0(x)\quad in~\Omega_f,\\
\mathcal{G}_f^2 [U_f^n,\Theta_f^n](x,t)&= 0\quad in~\Omega_f\times(0,T],\\
\Theta_f^n(x,t)&=0\quad on~\Gamma_M\times(0,T],\nonumber\\
\frac{\partial\Theta_f^n}{\partial n_f} &= 0 \quad on~\Gamma_E\times(0,T],\nonumber\\
\Theta_f^n(0,x) &= \Theta_f^0(x) \quad in~\Omega_f,\\
\mathcal{G}_p^1[U_p^n,P_p^n,\Theta_p^n](x,t)&=0\quad in~\Omega_p\times(0,T],\\
\nabla\cdot U_p^n &= 0 \quad in~\Omega_p\times(0,T],\nonumber\\
U_p^n(0,x) &= U_p^0(x) \quad in~\Omega_p,\nonumber\\
U_p^n\cdot n_p&=0\quad on~\partial\Omega_p\setminus\mathbb{I},\nonumber\\
\mathcal{G}_p^2[U_p^n,\Theta_p^n](x,t)&= 0 \quad in~\Omega_p\times(0,T],\\
\Theta_p^n(x,t)&= 0 \quad on~\Gamma_N\times(0,T],\nonumber\\
\frac{\partial\Theta_p^n}{\partial n_p} &= 0 \quad on~\Gamma_Z\times(0,T],\nonumber\\
\Theta_p^n(0,x) &= \Theta_p^0(x) \quad in~\Omega_p.
\end{align}

\noindent{{$\textbf{\emph{Theorem~2}}$}.~~}Under the assumptions of \emph{Hypotheses} $\textbf{\emph{H1}}$, $\textbf{\emph{Lemma1}}$ and $\textbf{\emph{Lemma2}}$, there exists a unique solution $(u_f,u_p;p_f,p_p;\theta_f,\theta_p)\in(Y_f\times Y_p\times Q_f\times Q_p\times W_T)$ to \eqref{equation-2.1}-\eqref{equation-2.16}. Moreover, when the sequence $(U_f^n,U_p^n,\Theta_f^n,\Theta_p^n)$ is uniformly bounded and equicontinuous, the neural networks $(U_f^n,U_p^n,\Theta_f^n,\Theta_p^n)$ converges strongly to $(u_f,u_p,\theta_f,\theta_p)$ in $L^2(\Omega_f)$ and $L^2(\Omega_p)$, respectively. Furthermore, $(U_f^n,U_p^n,\Theta_f^n,\Theta_p^n)$ uniformly converges to $(u_f, u_p,\theta_f,\theta_p)$ in $\Omega_{f}$ and $\Omega_{p}$, respectively.

The weak formulations of the coupled system \eqref{equation-2.1}-\eqref{equation-2.16} is to find $(U_f^n,U_p^n;P_f^n,P_p^n;\Theta_f^n,\Theta_p^n)\in [\mathfrak{C}_{u_f}^n(\varphi_f)]^d\times[\mathfrak{C}_{u_p}^n(\varphi_p)]^d\times\mathfrak{C}_p^n(\psi_f)\times\mathfrak{C}_p^n(\psi_p)\times\mathfrak{C}_{\theta}^n(\omega_f)\times\mathfrak{C}_{\theta}^n(\omega_p)$, for $\forall(V_f^n,V_p^n;Q_f^n,Q_p^n;W_f^n,W_p^n)\in  [\mathfrak{C}_{u_f}^n(\varphi_f)]^d\times[\mathfrak{C}_{u_p}^n(\varphi_p)]^d\times\mathfrak{C}_p^n(\psi_f)\times\mathfrak{C}_p^n(\psi_p)\times\mathfrak{C}_{\theta}^n(\omega_{\mathbb{I}})$ such that
\begin{align}
&(\frac{\partial U_f^n}{\partial t},V_f^n)_{\Omega_f}-(P_f^n, \nabla\cdot V_f^n)_{\Omega_f}+Pr(\nabla U_f^n,\nabla V_f^n)_{\Omega_f}+c(U_f^n, U_f^n, V_f^n)_{\Omega_f},\nonumber\\
&=PrRa\xi(\Theta_f^n,V_f^n)_{\Omega_f}+(f_f,V_f^n)_{\Omega_f},\label{equation-3.17}\\
&(\frac{\partial U_p^n}{\partial t},V_p^n)_{\Omega_p}+(U_p^n,V_p^n)_{\Omega_p}-Da(P_p^n,\nabla\cdot V_P^n)_{\Omega_p}=DaRa\xi(\Theta_p^n,V_p^n)_{\Omega_p},\\
&(Q_f^n,\nabla\cdot U_f^n)_{\Omega_f}=0,\\
&(Q_p^n,\nabla\cdot U_p^n)_{\Omega_p}=0,\\
&(\frac{\partial \Theta_f^n}{\partial t},W_f^n)_{\Omega_f}+k_f(\nabla\Theta_f^n,\nabla W_f^n)_{\Omega_f}+c_2(U_f^n, \Theta_f^n, W_f^n)_{\Omega_f}=(\Upsilon_f,W_f^n)_{\Omega_f},\\
&(\frac{\partial \Theta_p^n}{\partial t},W_p^n)_{\Omega_p}+k_p(\nabla\Theta_p^n,\nabla W_p^n)_{\Omega_p}+c_3(U_p^n, \Theta_p^n, W_p^n)_{\Omega_p}=(\Upsilon_p,W_p^n)_{\Omega_p}.\label{equation-3.22}
\end{align}

Let $\Theta_f^n=W_f^n=\Theta_f^n(0,x)$ in $(3.21)$, we get that
\begin{align}
&\Big(\frac{\partial \Theta_f^n(0,x)}{\partial t},\Theta_f^n(0,x)\Big)_{\Omega_f}+k_f\Big(\nabla\Theta_f^n(0,x),\nabla \Theta_f^n(0,x)\Big)_{\Omega_f}\nonumber\\
&+c_2\Big(U_f^n(0,x), \Theta_f^n(0,x), \Theta_f^n(0,x)\Big)_{\Omega_f}=\Big(\Upsilon_f(0,x),\Theta_f^n(0,x)\Big)_{\Omega_f}.
\end{align}
Using the Poincar$\acute{e}$ inequality and Young's inequality can lead to
\begin{align}
&\frac{d}{dt}\|\Theta_f^n(0,x)\|_{0(\Omega_f)}^2+c(k_f)\|\nabla \Theta_f^n(0,x)\|_{0(\Omega_f)}^2\leq \|\Upsilon_f(0,x)\|_{0(\Omega_f)}^2.
\end{align}
$\{\Theta_f^n(0,x)\}$ is bounded by (3.23) and (3.24), $\{\Theta_p^n(0,x)\},\{U_f^n(0,x)\}$ and $\{U_p^n(0,x)\}$ are bounded in the same way.

Let $\Theta_f^n=\Theta_f^n(0,x),W_f^n=\frac{\partial\Theta_f^n(0,x)}{\partial t}$ in $(3.21)$, we obtain that
\begin{align}
&\Big(\frac{\partial \Theta_f^n(0,x)}{\partial t},\frac{\partial\Theta_f^n(0,x)}{\partial t}\Big)_{\Omega_f}+k_f\Big(\nabla\Theta_f^n(0,x),\nabla \frac{\partial \Theta_f^n(0,x)}{\partial t}\Big)_{\Omega_f}\nonumber\\
&+c_2\Big(U_f^n(0,x), \Theta_f^n(0,x), \frac{\partial\Theta_f^n(0,x)}{\partial t}\Big)_{\Omega_f}=\Big(\Upsilon_f(0,x),\frac{\partial\Theta_f^n(0,x)}{\partial t}\Big)_{\Omega_f}.
\end{align}
Combine the Poincar$\acute{e}$, Cauchy-Schwarz and Young's inequalities, there is
\begin{align}
&\Big\|\frac{\partial\Theta_f^n(0,x)}{\partial t}\Big\|_{0(\Omega_f)}^2+c(k_f)\frac{d}{dt}\|\nabla \Theta_f^n(0,x)\|_{0(\Omega_f)}^2\leq \|\Upsilon_f(0,x)\|_{0(\Omega_f)}^2,
\end{align}
$\{\frac{\partial\Theta_f^n(0,x)}{\partial t}\}$ is bounded according to (3.25) and (3.26), so does $\{\frac{\partial\Theta_p^n(0,x)}{\partial t}\},\{\frac{\partial U_f^n(0,x)}{\partial t}\}$ and $\{\frac{\partial U_p^n(0,x)}{\partial t}\}$.
Add the equation (3.21) to (3.22), then let $W_f^n=\Theta_f^n$, $W_p^n=\Theta_p^n$, we have
\begin{align}
&\frac{d}{dt}(\|\Theta_f^n\|_{0(\Omega_f)}^2+\|\Theta_p^n\|_{0(\Omega_p)}^2)+(2k_f-c_1)\|\nabla \Theta_f^n\|_{0(\Omega_f)}^2+(2k_p-c_2)\|\nabla\Theta_p^n\|_{0(\Omega_p)}^2\nonumber\\
&\leq \|\Upsilon_f\|_{0(\Omega_f)}^2+\|\Upsilon_p\|_{0(\Omega_p)}^2,
\end{align}
and integrating the equation (3.27) in time, we get that
\begin{align}
&\|\Theta_f^n\|_{0(\Omega_f)}^2+\|\Theta_p^n\|_{0(\Omega_p)}^2+\int_0^t(2k_f-c_1)\|\nabla \Theta_f^n\|_{0(\Omega_f)}^2dl+\int_0^t(2k_p-c_2)\|\nabla\Theta_p^n\|_{0(\Omega_p)}^2dl\nonumber\\
&\leq \|\Theta_f^n(0,x)\|_{0(\Omega_f)}^2+\|\Theta_p^n(0,x)\|_{0(\Omega_p)}^2+\int_0^t\Big(\|\Upsilon_f\|_{0(\Omega_f)}^2+\|\Upsilon_p\|_{0(\Omega_p)}^2\Big)
\end{align}
This means that $\{\Theta_f^n\}$, $\{\Theta_p^n\}$ are uniformly bounded in $L^2(0,T;H^1(\Omega))$.

Plugging the equations (3.17) and (3.18), set $V_f^n=U_f^n$, $V_p^n=U_p^n$,
\begin{align}
&(\frac{\partial U_f^n}{\partial t},U_f^n)_{\Omega_f}+(\frac{\partial U_p^n}{\partial t},U_p^n)_{\Omega_p}-(P_f^n, \nabla\cdot U_f^n)_{\Omega_f}+Pr(\nabla U_f^n,\nabla U_f^n)_{\Omega_f}\nonumber\\
&+c_1(U_f^n, U_f^n, U_f^n)_{\Omega_f}+(U_p^n,U_p^n)_{\Omega_p}-Da(P_p^n,\nabla\cdot U_P^n)_{\Omega_p}\nonumber\\
&=PrRa\xi(\Theta_f^n,U_f^n)_{\Omega_f}+DaRa\xi(\Theta_p^n,U_p^n)_{\Omega_p}+(f_f,U_f^n)_{\Omega_f}.
\end{align}
Integrating the equation (3.29) in time and using the Gronwall inequality can lead to
\begin{align}
&\|U_f^n\|_{0(\Omega_f)}^2+\|U_p^n\|_{0(\Omega_p)}^2+2(Pr-c_3)\int_0^t\|\nabla U_f^n\|_{0(\Omega_f)}^2dl+2\int_0^t\|U_p^n\|_{0(\Omega_p)}^2dl\nonumber\\
\leq&
c(Pr,Ra,Da)\int_0^t\Big(\|U_f^n\|_{0(\Omega_f)}^2+\|U_p^n\|_{0(\Omega_p)}^2+\|\Theta_f^n\|_{0(\Omega_f)}^2+\|\Theta_p^n\|_{0(\Omega_p)}^2\Big)dl\nonumber\\
&+ \|U_f^n(0,x)\|_{0(\Omega_f)}^2+\|U_p^n(0,x)\|_{0(\Omega_p)}^2+\int_0^t\|f_f\|_{0(\Omega_f)}^2dl\nonumber\\
\leq& \Big(\|U_f^n(0,x)\|_{0(\Omega_f)}^2+M_1(t)\Big)exp\Big(\int_0^tc(Pr,Ra,Da)dl\Big),
\end{align}
where,
$M_1(t)=\int_0^t\|f_f\|_{0(\Omega_f)}^2+c(Pr,Ra,Da)\int_0^t\Big(\|\Theta_f^n\|_{0(\Omega_f)}^2+\|\Theta_p^n\|_{0(\Omega_p)}^2\Big).$\\
This shows that $\{U_f^n\}$ and $\{U_p^n\}$ are uniformly bounded in $L^2(0,T;H^1(\Omega))$ and $L^2(0,T;L^2(\Omega))$ respectively.
Add (3.17) to (3.21) and let $V_f^n=U_f^n, W_f^n=\Theta_f^n$ we have
\begin{align}
&\frac{d}{dt}\Big(\Big\|\frac{\partial U_f^n}{\partial t}\Big\|_{L^2(\Omega_f)}^2+\Big\|\frac{\partial \Theta_f^n}{\partial t}\Big\|_{L^2(\Omega_f)}^2\Big)+2\Big(Pr\Big\|\nabla\frac{\partial U_f^n}{\partial t} \Big\|_{L^2(\Omega_f)}^2+k_f\Big\|\nabla\frac{\partial \Theta_f^n}{\partial t}\Big\|_{L^2(\Omega_f)}^2\Big)\nonumber\\
\leq& \Big|2\int_{\Omega_f}\Big(\frac{\partial U_f^n}{\partial t}\cdot\nabla\Big)U_f^n\cdot\frac{\partial U_f^n}{\partial t}d\Omega_f\Big|+\Big|2\int_{\Omega_f}\Big(U_f^n\cdot\nabla\Big)\frac{\partial U_f^n}{\partial t}\cdot\frac{\partial U_f^n}{\partial t}d\Omega_f\Big|\nonumber\\
&+\Big|2\int_{\Omega_f}\frac{\partial U_f^n}{\partial t}\cdot\nabla\Theta_f^n\cdot\frac{\partial\Theta_f^n}{\partial t}d\Omega_f\Big|+\Big|2\int_{\Omega_f}U_f^n\cdot\nabla\frac{\partial\Theta_f^n}{\partial t}\cdot\frac{\partial\Theta_f^n}{\partial t}d\Omega_f\Big|\nonumber\\
&+\Big|2PrRa\xi\int_{\Omega_f}\frac{\partial\Theta_f^n}{\partial t}\cdot\frac{\partial U_f^n}{\partial t}d\Omega_f\Big|+\Big|2\int_{\Omega_f}\frac{\partial f_f}{\partial t}\cdot\frac{\partial U_f^n}{\partial t}d\Omega_f\Big|+\Big|2\int_{\Omega_f}\frac{\partial\Upsilon_f}{\partial t}\cdot\frac{\partial \Theta_f^n}{\partial t}d\Omega_f\Big|.
\end{align}
Moreover, we will estimate the right hand side terms of (3.31),
\begin{align}
\Big|2\int_{\Omega_f}\Big(\frac{\partial U_f^n}{\partial t}\cdot\nabla\Big)U_f^n\cdot\frac{\partial U_f^n}{\partial t}d\Omega_f\Big|\leq C_1(Pr)\Big\|\frac{\partial U_f^n}{\partial t}\Big\|_{L^2(\Omega_f)}^2\Big\|\nabla U_f^n\Big\|_{L^2(\Omega_f)}^2+\frac{Pr}{2}\Big\|\frac{\partial U_f^n}{\partial t}\Big\|_{L^2(\Omega_f)}^2,
\end{align}

\begin{align}
\Big|2\int_{\Omega_f}\Big(U_f^n\cdot\nabla\Big)\frac{\partial U_f^n}{\partial t}\cdot\frac{\partial U_f^n}{\partial t}d\Omega_f\Big|\leq C_2(Pr)\Big\|U_f^n\Big\|_{L^2(\Omega_f)}^2\Big\|\nabla\frac{\partial U_f^n}{\partial t}\Big\|_{L^2(\Omega_f)}^2+\frac{Pr}{2}\Big\|\frac{\partial U_f^n}{\partial t}\Big\|_{L^2(\Omega_f)}^2,
\end{align}

\begin{align}
\Big|2\int_{\Omega_f}\frac{\partial U_f^n}{\partial t}\cdot\nabla\Theta_f^n\cdot\frac{\partial\Theta_f^n}{\partial t}d\Omega_f\Big|\leq& C_1(Pr,k_f)\Big(\Big\|\frac{\partial U_f^n}{\partial t}\Big\|_{L^2(\Omega_f)}^2\Big\|\nabla \Theta_f^n\Big\|_{L^2(\Omega_f)}^2+\Big\|\frac{\partial \Theta_f^n}{\partial t}\Big\|_{L^2(\Omega_f)}^2\Big\|\nabla \Theta_f^n\Big\|_{L^2(\Omega_f)}^2\Big)\nonumber\\
&+\Big(\frac{Pr}{2}\Big\|\nabla\frac{\partial U_f^n}{\partial t}\Big\|_{L^2(\Omega_f)}^2+\frac{k_f}{2}\Big\|\nabla\frac{\partial \Theta_f^n}{\partial t}\Big\|_{L^2(\Omega_f)}^2\Big),
\end{align}

\begin{align}
\Big|2\int_{\Omega_f}U_f^n\cdot\nabla\frac{\partial\Theta_f^n}{\partial t}\cdot\frac{\partial\Theta_f^n}{\partial t}d\Omega_f\Big|\leq& C_2(Pr,k_f)\Big(\Big\|U_f^n\Big\|_{L^2(\Omega_f)}^2\Big\|\nabla\frac{\partial\Theta_f^n}{\partial t}\Big\|_{L^2(\Omega_f)}^2+\Big\|\frac{\partial \Theta_f^n}{\partial t}\Big\|_{L^2(\Omega_f)}^2\Big\|\nabla\frac{\partial\Theta_f^n}{\partial t}\Big\|_{L^2(\Omega_f)}^2\Big)\nonumber\\
&+\Big(\frac{Pr}{2}\Big\|\nabla U_f^n\Big\|_{L^2(\Omega_f)}^2+\frac{k_f}{2}\Big\|\nabla\frac{\partial \Theta_f^n}{\partial t}\Big\|_{L^2(\Omega_f)}^2\Big),
\end{align}
and
\begin{align}
\Big|2PrRa\xi\int_{\Omega_f}\frac{\partial\Theta_f^n}{\partial t}\cdot\frac{\partial U_f^n}{\partial t}d\Omega_f\Big|\leq C(Pr,Ra)\Big(\Big\|\frac{\partial\Theta_f^n}{\partial t}\Big\|_{L^2(\Omega_f)}^2+\Big\|\frac{\partial U_f^n}{\partial t}\Big\|_{L^2(\Omega_f)}^2\Big).
\end{align}
Plugging the above estimates into (3.31), we obtain that
\begin{align}
&\frac{d}{dt}\Big(\Big\|\frac{\partial U_f^n}{\partial t}\Big\|_{L^2(\Omega_f)}^2+\Big\|\frac{\partial \Theta_f^n}{\partial t}\Big\|_{L^2(\Omega_f)}^2\Big)+2\Big(Pr\Big\|\nabla\frac{\partial U_f^n}{\partial t} \Big\|_{L^2(\Omega_f)}^2+k_f\Big\|\nabla\frac{\partial \Theta_f^n}{\partial t}\Big\|_{L^2(\Omega_f)}^2\Big)\nonumber\\
\leq& C_1(Pr)\Big\|\frac{\partial U_f^n}{\partial t}\Big\|_{L^2(\Omega_f)}^2\Big\|\nabla U_f^n\Big\|_{L^2(\Omega_f)}^2+\frac{Pr}{2}\Big\|\frac{\partial U_f^n}{\partial t}\Big\|_{L^2(\Omega_f)}^2\nonumber\\
&+C_2(Pr)\Big\|U_f^n\Big\|_{L^2(\Omega_f)}^2\Big\|\nabla\frac{\partial U_f^n}{\partial t}\Big\|_{L^2(\Omega_f)}^2+\frac{Pr}{2}\Big\|\frac{\partial U_f^n}{\partial t}\Big\|_{L^2(\Omega_f)}^2\nonumber\\
&+C_1(Pr,k_f)\Big(\Big\|\frac{\partial U_f^n}{\partial t}\Big\|_{L^2(\Omega_f)}^2\Big\|\nabla \Theta_f^n\Big\|_{L^2(\Omega_f)}^2+\Big\|\frac{\partial \Theta_f^n}{\partial t}\Big\|_{L^2(\Omega_f)}^2\Big\|\nabla \Theta_f^n\Big\|_{L^2(\Omega_f)}^2\Big)\nonumber\\
&+\Big(\frac{Pr}{2}\Big\|\nabla\frac{\partial U_f^n}{\partial t}\Big\|_{L^2(\Omega_f)}^2+\frac{k_f}{2}\Big\|\nabla\frac{\partial \Theta_f^n}{\partial t}\Big\|_{L^2(\Omega_f)}^2\Big)\nonumber\\
&+C_2(Pr,k_f)\Big(\Big\|U_f^n\Big\|_{L^2(\Omega_f)}^2\Big\|\nabla\frac{\partial\Theta_f^n}{\partial t}\Big\|_{L^2(\Omega_f)}^2+\Big\|\frac{\partial \Theta_f^n}{\partial t}\Big\|_{L^2(\Omega_f)}^2\Big\|\nabla\frac{\partial\Theta_f^n}{\partial t}\Big\|_{L^2(\Omega_f)}^2\Big)\nonumber\\
&+\Big(\frac{Pr}{2}\Big\|\nabla U_f^n\Big\|_{L^2(\Omega_f)}^2+\frac{k_f}{2}\Big\|\nabla\frac{\partial \Theta_f^n}{\partial t}\Big\|_{L^2(\Omega_f)}^2\Big)+C(Pr,Ra)\Big(\Big\|\frac{\partial\Theta_f^n}{\partial t}\Big\|_{L^2(\Omega_f)}^2+\Big\|\frac{\partial U_f^n}{\partial t}\Big\|_{L^2(\Omega_f)}^2\Big)\nonumber\\
&+\Big\|\frac{\partial f_f}{\partial t}\Big\|_{L^2(\Omega_f)}^2+\Big\|\frac{\partial U_f^n}{\partial t}\Big\|_{L^2(\Omega_f)}^2+\Big\|\frac{\partial \Upsilon_f}{\partial t}\Big\|_{L^2(\Omega_f)}^2+\Big\|\frac{\partial \Theta_f^n}{\partial t}\Big\|_{L^2(\Omega_f)}^2.
\end{align}
Rearranging the equation (3.37), integrating the equation in time and using the Gronwall inequality, we can get that
\begin{align}
&\Big\|\frac{\partial U_f^n}{\partial t}\Big\|_{L^2(\Omega_f)}^2+\Big\|\frac{\partial \Theta_f^n}{\partial t}\Big\|_{L^2(\Omega_f)}^2\nonumber\\
\leq& \Big\|\frac{\partial U_f^n(0,x)}{\partial t}\Big\|_{L^2(\Omega_f)}^2+\Big\|\frac{\partial \Theta_f^n(0,x)}{\partial t}\Big\|_{L^2(\Omega_f)}^2+M_1(t)+C_2(Pr,k_f)\Big\|U_f^n\Big\|_{L^2(\Omega_f)}^2\int_0^t\Big\|\nabla\frac{\partial\Theta_f^n}{\partial t}\Big\|_{L^2(\Omega_f)}^2ds\nonumber\\
&+\Big(C_1(Pr)\Big\|\nabla U_f^n\Big\|_{L^2(\Omega_f)}^2+C_2(Pr)\|U_f^n\|_{L^2}^2+C_1(Pr,k_f)\Big\|\nabla \Theta_f^n\Big\|_{L^2(\Omega_f)}^2+C(Pr,Ra)\Big)\int_0^t\Big\|\nabla\frac{\partial U_f^n}{\partial t}\Big\|_{L^2(\Omega_f)}^2ds\nonumber\\
&+\Big(C_2(Pr,k_f)+C(Pr,Ra)\Big)\int_0^t\Big\|\nabla\frac{\partial\Theta_f^n}{\partial t}\Big\|_{L^2(\Omega_f)}^2\Big\|\frac{\partial \Theta_f^n}{\partial t}\Big\|_{L^2(\Omega_f)}^2ds\nonumber\\
\leq&\Big(\Big\|\frac{\partial U_f^n(0,x)}{\partial t}\Big\|_{L^2(\Omega_f)}^2+\Big\|\frac{\partial \Theta_f^n(0,x)}{\partial t}\Big\|_{L^2(\Omega_f)}^2+M_1(t)\Big)exp({\int_0^t\tilde{C})ds},
\end{align}
where
$$M_1(t)=\int_0^t\Big(\frac{Pr}{2}\Big\|\nabla U_f^n\Big\|_{L^2(\Omega_f)}^2+\Big\|\frac{\partial f_f}{\partial t}\Big\|_{L^2(\Omega_f)}^2+\Big\|\frac{\partial \Upsilon_f}{\partial t}\Big\|_{L^2(\Omega_f)}^2\Big)ds,$$
$$\tilde{C}=C_1(Pr)\Big\|\nabla U_f^n\Big\|_{L^2(\Omega_f)}^2+C_2(Pr)\|U_f^n\|_{L^2}^2+C_1(Pr,k_f)\Big\|\nabla \Theta_f^n\Big\|_{L^2(\Omega_f)}^2+C(Pr,Ra),$$
This shows that $\Big\{\frac{\partial U_f^n}{\partial t}\Big\}$ and $\Big\{\frac{\partial \Theta_f^n}{\partial t}\Big\}$ are uniformly bounded in $L^2\Big(0,T;H^1(\Omega_f)\Big)$.
Add (3.18) to (3.22) and let $V_p^n=U_p^n, W_p^n=\Theta_p^n$, we have
\begin{align}
&\frac{d}{dt}\Big(\Big\|\frac{\partial U_p^n}{\partial t}\Big\|_{L^2(\Omega_p)}^2+\Big\|\frac{\partial\Theta_p^n}{\partial t}\Big\|_{L^2(\Omega_p)}^2\Big)+2k_p\Big\|\nabla\frac{\partial \Theta_p^n}{\partial t}\Big\|_{L^2(\Omega_p)}^2+2\Big\|\frac{\partial U_p^n}{\partial t}\Big\|_{L^2(\Omega_p)}^2\nonumber\\
\leq& \Big|2DaRa\xi\int_{\Omega_p}\frac{\partial\Theta_p^n}{\partial t}\cdot\frac{\partial U_p^n}{\partial t}d\Omega_p\Big|+\Big|2\int_{\Omega_p}\frac{\partial U_p^n}{\partial t}\cdot\nabla\Theta_p^n\cdot\frac{\partial\Theta_p^n}{\partial t}d\Omega_p\Big|\nonumber\\
&+\Big|2\int_{\Omega_p}U_p^n\cdot\nabla\frac{\partial\Theta_p^n}{\partial t}\cdot\frac{\partial\Theta_p^n}{\partial t}d\Omega_p\Big|+\Big|2\int_{\Omega_p}\frac{\partial\Upsilon_p}{\partial t}\cdot\frac{\partial \Theta_p^n}{\partial t}d\Omega_p\Big|.
\end{align}
In the following we will estimate the terms on the right hand side of (3.39) one by one as follows.
\begin{align}
\Big|2\int_{\Omega_p}\frac{\partial U_p^n}{\partial t}\cdot\nabla\Theta_p^n\cdot\frac{\partial\Theta_p^n}{\partial t}d\Omega_p\Big|\leq C(k_p)\Big\|\frac{\partial U_p^n}{\partial t}\|_{L^2(\Omega_p)}^2\Big\|\Theta_p^n\Big\|_{L^2(\Omega_p)}^2+\frac{k_p}{2}\Big\|\nabla\frac{\partial\Theta_p^n}{\partial t}\Big\|_{L^2(\Omega_p)}^2,
\end{align}

\begin{align}
\Big|2\int_{\Omega_p}U_p^n\cdot\nabla\frac{\partial\Theta_p^n}{\partial t}\cdot\frac{\partial\Theta_p^n}{\partial t}d\Omega_p\Big|\leq C(k_p)\Big\|U_p^n\Big\|_{L^2(\Omega_p)}^2\Big\|\frac{\partial\Theta_p^n}{\partial t}\Big\|_{L^2(\Omega_p)}^2+\frac{k_p}{2}\Big\|\nabla\frac{\partial\Theta_p^n}{\partial t}\Big\|_{L^2(\Omega_p)}^2,
\end{align}

\begin{align}
\Big|2DaRa\xi\int_{\Omega_p}\frac{\partial\Theta_p^n}{\partial t}\cdot\frac{\partial U_p^n}{\partial t}d\Omega_p\Big|\leq C(Da,Ra)\Big\|\frac{\partial\Theta_p^n}{\partial t}\Big\|_{L^2(\Omega_p)}^2+\Big\|\frac{\partial U_p^n}{\partial t}\Big\|_{L^2(\Omega_p)}^2,
\end{align}

\begin{align}
\Big|2\int_{\Omega_p}\frac{\partial\Upsilon_p}{\partial t}\cdot\frac{\partial \Theta_p^n}{\partial t}d\Omega_p\Big|\leq C\Big(\Big\|\frac{\partial\Upsilon_p}{\partial t}\Big\|_{L^2(\Omega_p)}^2+\Big\|\frac{\partial \Theta_p^n}{\partial t}\Big\|_{L^2(\Omega_p)}^2\Big),
\end{align}
submitting $(3.40)-(3.43)$ to $(3.29)$, we obtain that
\begin{align}
&\frac{d}{dt}\Big(\Big\|\frac{\partial U_p^n}{\partial t}\Big\|_{L^2(\Omega_p)}^2+\Big\|\frac{\partial\Theta_p^n}{\partial t}\Big\|_{L^2(\Omega_p)}^2\Big)+k_p\Big\|\nabla\frac{\partial \Theta_p^n}{\partial t}\Big\|_{L^2(\Omega_p)}^2\nonumber\\
\leq& C\Big\|\frac{\partial\Upsilon_p}{\partial t}\Big\|_{L^2(\Omega_p)}^2+\Big(C(k_p)\Big\|\Theta_p^n\Big\|_{L^2(\Omega_p)}^2-1\Big)\Big\|\frac{\partial U_p^n}{\partial t}\Big\|_{L^2(\Omega_p)}^2\nonumber\\
&+\Big(C(k_p)+C(Da,Ra)+C\Big)\Big\|\frac{\partial \Theta_p^n}{\partial t}\Big\|_{L^2(\Omega_p)}^2.
\end{align}
Integrating the equation (3.44) in time and using the Gronwall inequality can lead to
\begin{align}
&\Big\|\frac{\partial U_p^n}{\partial t}\Big\|_{L^2(\Omega_p)}^2+\Big\|\frac{\partial\Theta_p^n}{\partial t}\Big\|_{L^2(\Omega_p)}^2+k_p\int_0^t\Big\|\nabla\frac{\partial \Theta_p^n}{\partial t}\Big\|_{L^2(\Omega_p)}^2dl\nonumber\\
\leq& \Big\|\frac{\partial U(0,x)_p^n}{\partial t}\Big\|_{L^2(\Omega_p)}^2+\Big\|\frac{\partial\Theta(0,x)_p^n}{\partial t}\Big\|_{L^2(\Omega_p)}^2+C\int_0^t\Big\|\frac{\partial\Upsilon_p}{\partial t}\Big\|_{L^2(\Omega_p)}^2dl\nonumber\\
&+\int_0^tC(k_p)\Big\|\Theta_p^n\Big\|_{L^2(\Omega_p)}^2\Big\|\frac{\partial U_p^n}{\partial t}\Big\|_{L^2(\Omega_p)}^2dl+\int_0^t\Big(C(k_p)+C(Da,Ra)+C\Big)\Big\|\frac{\partial \Theta_p^n}{\partial t}\Big\|_{L^2(\Omega_p)}^2dl\nonumber\\
\leq& \Big(\Big\|\frac{\partial U(0,x)_p^n}{\partial t}\Big\|_{L^2(\Omega_p)}^2+\Big\|\frac{\partial\Theta(0,x)_p^n}{\partial t}\Big\|_{L^2(\Omega_p)}^2+M_2(t)\Big)exp({\int_0^tC(k_p,Da,Ra)dl}),
\end{align}
where
$$M_2(t)=C\int_0^t\Big\|\frac{\partial\Upsilon_p}{\partial t}\Big\|_{L^2(\Omega_p)}^2dl.$$
This shows that $\Big\{\frac{\partial U_p^n}{\partial t}\Big\}$ and $\Big\{\frac{\partial \theta_p^n}{\partial t}\Big\}$ are uniformly bounded in $L^2\Big(0,T;H^2(\Omega_p)\Big)$.

From the above proof, we can get that $\{U_f^n\}$ and $\Big\{\frac{\partial U_f^n}{\partial t}\Big\}$ are uniformly bounded in $L^2(0,T;H^1(\Omega_f))$, combining with the Aubin-Lions's compactness lemma, there exists a subsequence of $\{U_f^n\}$ (still denoted by the $\{U_f^n\}$), which converges to $u_f\in L^{\infty}(0,T;L^2(\Omega_f))\cap L^2(0,T;H_0^1(\Omega_f))$ such that
$$U_f^n\rightarrow u_f~in~L^2(0,T;L^2(\Omega_f)).$$
We can get the
$$\Theta_f^n\rightarrow \theta_f~in~L^2(0,T;L^2(\Omega_f)),$$
$$U_p^n\rightarrow u_p~in~L^2(0,T;H^{-1}(\Omega_p)),$$
$$\Theta_p^n\rightarrow~\theta_p~in~L^2(0,T;H^{-1}(\Omega_p)).$$

Up to now, we have been ready for passing to the limit as $n\rightarrow\infty$ in the weak formulation, it is easy to see that $u_f,\theta_f,u_p$ and $\theta_p$ satisfy \eqref{equation-3.17}-\eqref{equation-3.22} in a weak sense. This result shows that $U_f^n$ and $\Theta_f^n$ converge strongly to $u_f$ and $\theta_f$ in $L^2(\Omega_f)$, respectively. $U_p^n$ and $\Theta_p^n$ converge strongly to $u_p$ and $\theta_p$ in $H^{-1}(\Omega_p)$, respectively.  Furthermore, based on the equicontinuity and uniform boundedness of $U_f^n$ and $\Theta_f^n$, we obtain that $U_f^n$ and $\Theta_f^n$ uniformly converges to $u_f$ and $\theta_f$ in $\Omega_f$ by Arzel$\grave{a}$-Ascoli theorem. Similarly, it can be obtained that $U_p^n$ and $\Theta_p^n$ uniformly converges to $u_p$ and $\theta_p$ in $\Omega_p$.

\section{Numerical results}

In this section, we present two numerical experiments to validate the proposed DGM for the closed-loop geothermal system. In section 4.1, an analytic solution is provided to show the convergence of the neural network solutions and the impact of the hidden layers on the convergence. Section 4.2 is conducted to investigate a benchmark problem for thermal convection in a squared cavity.

\subsection{Convergence and stability tests.}

In this examples, we consider the problem \eqref{equation-2.1}-\eqref{equation-2.16} on a $2D$ spatial domain consisted of two subdomains. One is the free fluid flow region $\Omega_f=[0,1]\times[1,2]$, and the other one is the porous media region $\Omega_p = [0,1] \times[0,1]$. The selected exact solution for the model is given by

\begin{align}
u_f &= \begin{pmatrix} 10x^2(x -1)^2y(y - 1)(2y - 1)cos(t)  \\ -10x(x - 1)(2x - 1)y^2(y - 1)^2cos(t) \end{pmatrix},\nonumber\\
p_f &= 10(2x - 1)(2y - 1)cos(t),\nonumber\\
u_p & = \begin{pmatrix} [2\pi sin2(\pi x)sin(\pi y)cos(\pi y)]cos(t)  \\ [-2\pi sin(\pi x)sin2(\pi y)cos(\pi x)]cos(t) \end{pmatrix},\nonumber\\
\phi_p &=cos(\pi x)cos(\pi y)cos(t),\nonumber\\
\theta_f &= ax(1-x)(1-y)e^{-t},\nonumber\\
\theta_p &= ax(1-x)(y-y^2)e^{-t}.
\end{align}

The initial condition, boundary condition, and forcing term of the model problem \eqref{equation-2.1} to \eqref{equation-2.16} are determined
by the exact solutions of the model problem. Choose the parameters value $a=1.0$, $Pr=1.0$, $Ra=1.0$, $k_f=k_p=1.0$, $Da=1.0$ and  $\gamma=10^5$. For convenience, we denote the $L^2$ error ($errL^2$) and relative $L^2$ error ($errRL^2$) as follows
$$errL^2=\frac{1}{m}\sum\limits_{i=1}^m(U_i-u_i)^2,\quad errL^2=\frac{1}{m}\sum\limits_{i=1}^m(\Theta_i-\theta_i)^2,$$
$$errRL^2=\frac{\frac{1}{m}\sum\limits_{i=1}^m(U_i-u_i)^2}{\frac{1}{m}\sum\limits_{i=1}^m(u_i)^2},\quad errRL^2=\frac{\frac{1}{m}\sum\limits_{i=1}^m(\Theta_i-\theta_i)^2}{\frac{1}{m}\sum\limits_{i=1}^m(\theta_i)^2},$$
where $U_i$ and $\Theta_i$($i=1,2,\cdots,m$) are the neural network solutions of velocity and temperature, respectively. And $u_i$ and $\theta_i$($i=1,2,\cdots,m$) are the corresponding exact solutions of velocity and temperature.

\begin{figure}[H]
  \centering
  \subfigure{
  \begin{minipage}[t]{0.4\linewidth}
  \centering
  \includegraphics[width=2.7in]{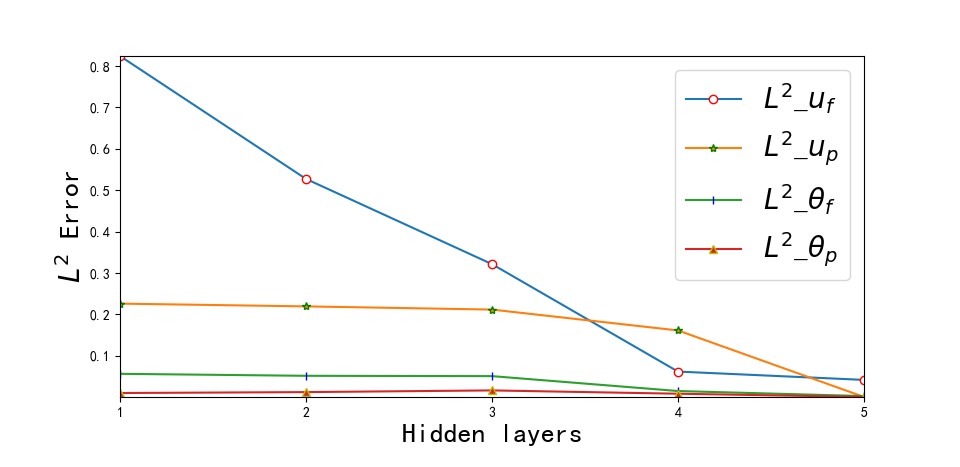}\\
  \end{minipage}%
}
  \subfigure{
  \begin{minipage}[t]{0.4\linewidth}
  \centering
  \includegraphics[width=2.7in]{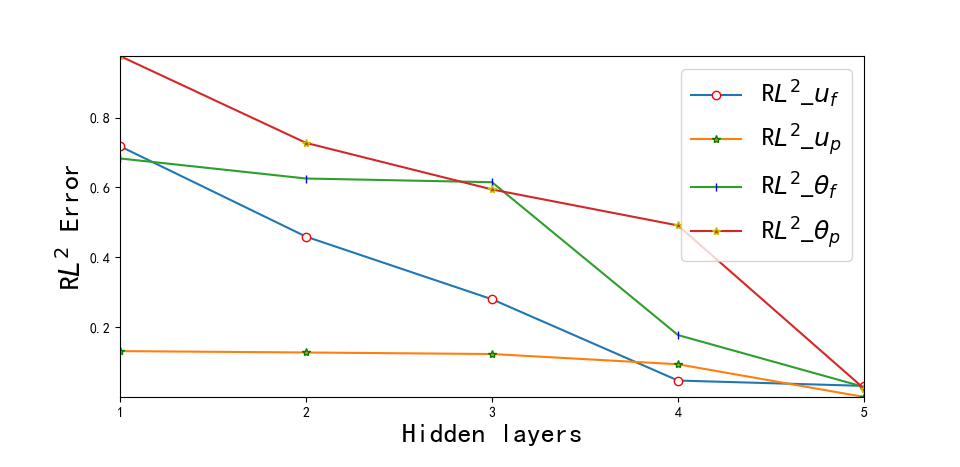}\\
  \end{minipage}%
}
\caption{Left: the $L^2$ error of different hidden layers; Right: the relative $L^2$ error of different hidden layers.}\label{figure-figure3}
\end{figure}

\begin{figure}[H]
  \centering
  \subfigure[exact solution]{
  \begin{minipage}[t]{0.2\linewidth}
  \centering
  \includegraphics[width=1.4in]{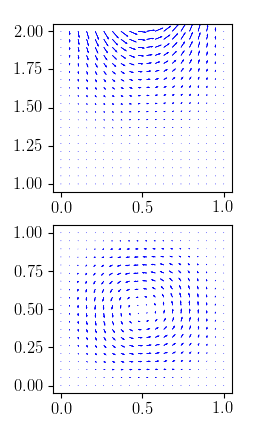}
  \end{minipage}%
}
  \subfigure[1 hidden layers]{
  \begin{minipage}[t]{0.2\linewidth}
  \centering
  \includegraphics[width=1.4in]{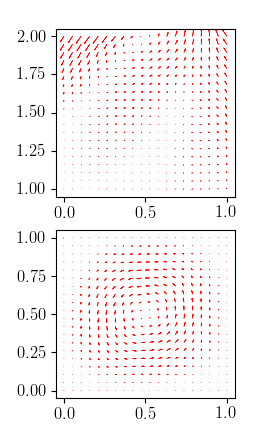}
  \end{minipage}%
}
  \subfigure[2 hidden layers]{
  \begin{minipage}[t]{0.2\linewidth}
  \centering
  \includegraphics[width=1.4in]{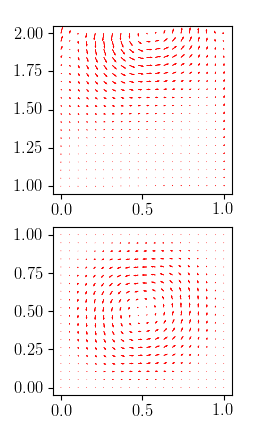}
  \end{minipage}%
}

  \subfigure[3 hidden layers]{
  \begin{minipage}[t]{0.2\linewidth}
  \centering
  \includegraphics[width=1.4in]{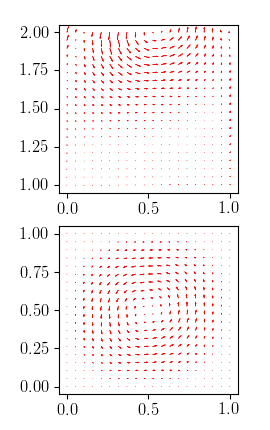}
  \end{minipage}%
}
  \subfigure[4 hidden layers]{
  \begin{minipage}[t]{0.2\linewidth}
  \centering
  \includegraphics[width=1.4in]{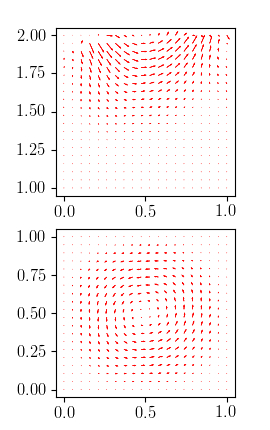}
  \end{minipage}%
}
  \subfigure[5 hidden layers]{
  \begin{minipage}[t]{0.2\linewidth}
  \centering
  \includegraphics[width=1.4in]{u_4_16.png}
  \end{minipage}%
}
\centering
\caption{A comparison between exact solution and neural network solution of velocity. }\label{figure-figure4}
\end{figure}

\begin{figure}[H]
  \centering
  \subfigure[exact solution]{
  \begin{minipage}[t]{0.2\linewidth}
  \centering
  \includegraphics[width=1.3in]{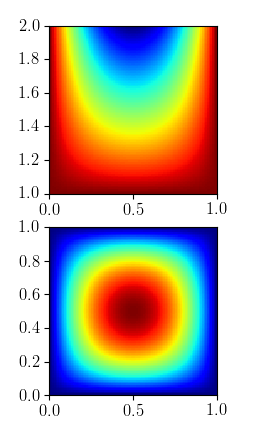}
  \end{minipage}%
}
  \subfigure[1 hidden layers]{
  \begin{minipage}[t]{0.2\linewidth}
  \centering
  \includegraphics[width=1.3in]{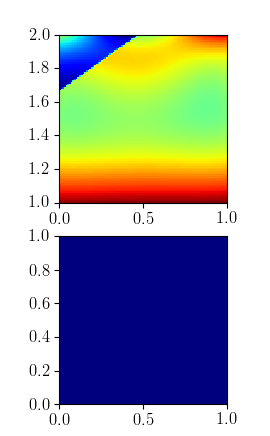}
  \end{minipage}%
}
  \subfigure[2 hidden layers]{
  \begin{minipage}[t]{0.2\linewidth}
  \centering
  \includegraphics[width=1.3in]{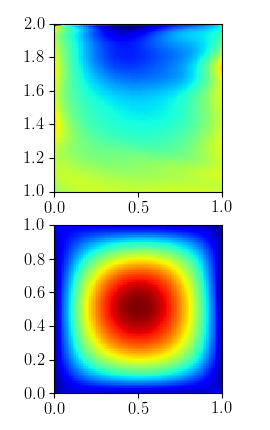}
  \end{minipage}%
}

  \subfigure[3 hidden layers]{
  \begin{minipage}[t]{0.2\linewidth}
  \centering
  \includegraphics[width=1.3in]{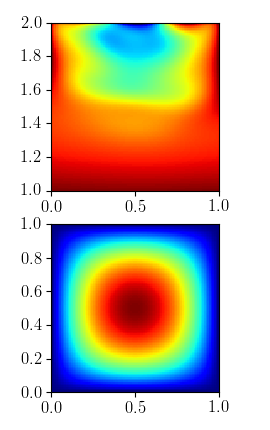}
  \end{minipage}%
}
  \subfigure[4 hidden layers]{
  \begin{minipage}[t]{0.2\linewidth}
  \centering
  \includegraphics[width=1.3in]{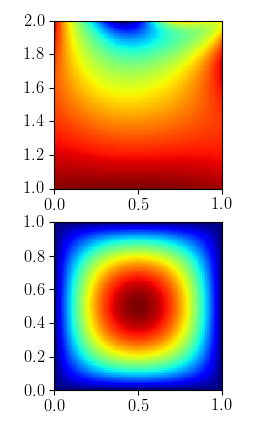}
  \end{minipage}%
}
  \subfigure[5 hidden layers]{
  \begin{minipage}[t]{0.2\linewidth}
  \centering
  \includegraphics[width=1.3in]{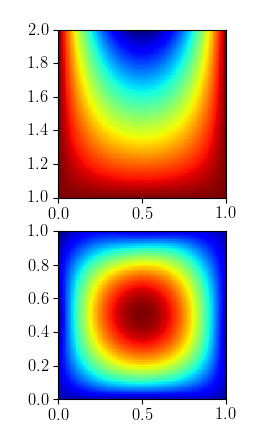}
  \end{minipage}%
}
\centering
\caption{A comparison between exact solution and neural network solution of temperature.}\label{figure-figure5}
\end{figure}

We train five architectures of the neural networks on the training set with 900 sample points, these neural networks include one to five hidden layers respectively and each hidden layer has 16 units. Table \ref{table-table1} and Figure \ref{figure-figure3} demonstrate the approximation ability of the DGM for the closed-loop geothermal system \eqref{equation-2.1}-\eqref{equation-2.16} on $2D$ spatial domain. Moreover, the value of the loss function $J(\overline{U})$ decreases gradually while the number of hidden layers increases, which is consistent with the result of theorem 1. Besides, combining the $L^2$ error and relative $L^2$ error, we find that when the loss function goes to zero, the solution of the neural network converges to the exact solution, which is the content expressed in theorem 2. Obviously, Figure \ref{figure-figure4} and Figure \ref{figure-figure5} show the good agreement between the exact solution and neural networks results with the increasing of the hidden layers. These results suggest that the neural network is close to the exact solution with the increasing of hidden layers, which is consistent with the theory.

\begin{table}[H]
  \caption{Numerical results of neural network. }
  \centering
  \subtable[$errL^2$]{
  \begin{tabular}{cc cc cc cc }
  \hline
  hidden layers & 1 & 2 & 3 & 4 & 5 \\
  \hline
  Dataset & 900 & 900 & 900 & 900 & 900 \\
  \hline
  $u_f$ & 0.8243 & 0.5271 & 0.3212 & 0.0616 & 0.0415 \\
  $u_p$ & 0.2258 & 0.2191 & 0.2114 & 0.1610& 0.0003 \\
  $\theta_f$ & 0.0562 & 0.0514 & 0.0506 & 0.0146 & 0.0024 \\
  $\theta_p$ & 0.0161 & 0.0120 & 0.0098 & 0.0081 & 0.0004 \\
  \hline
  \end{tabular}
}
  \qquad

  \subtable[$errRL^2$]{
  \begin{tabular}{cc cc cc cc }
  \hline
  hidden layers & 1 & 2 & 3 & 4 & 5 \\
  \hline
  Dataset & 900 & 900 & 900 & 900 & 900 \\
  \hline
  $u_f$ & 0.7177 & 0.4590 & 0.2797 & 0.0470 & 0.0317 \\
  $u_p$ & 0.1313 & 0.1273 & 0.1229 & 0.0936 & 0.00002 \\
  $\theta_f$ & 0.6828 & 0.6245 & 0.6148 & 0.1772 & 0.0287 \\
  $\theta_p$ & 0.9760 & 0.7274 & 0.5941 & 0.4910 & 0.0242 \\
  \hline
  \end{tabular}
}
  \qquad
  \subtable[$J(\overline{U})$]{
  \begin{tabular}{cc cc cc cc }
  \hline
  hidden layers & 1 & 2 & 3 & 4 & 5 \\
  \hline
  Dataset & 900 & 900 & 900 & 900 & 900 \\
  \hline
  $J_f$ & 0.6651 & 0.4510 & 0.2346 & 0.1811 & 0.0587\\
  $J_p$ & 0.1984 & 0.1406 & 0.0972 & 0.0387 & 0.0201\\
  $J_{\mathbb{I}}$ & 0.0127 & 0.0038 & 0.0040 & 0.0017 & 0.0013\\
  $J(\overline{U})$ & 0.8762 & 0.6043 & 0.3356 & 0.2215 & 0.0801\\
  \hline
  \end{tabular}
 }\label{table-table1}
\end{table}

\subsection{Convection in a squared cavity.}

\begin{figure}[H]
  \centering
  \subfigure[]{
  \begin{minipage}[t]{0.3\linewidth}
  \centering
  \includegraphics[width=2in]{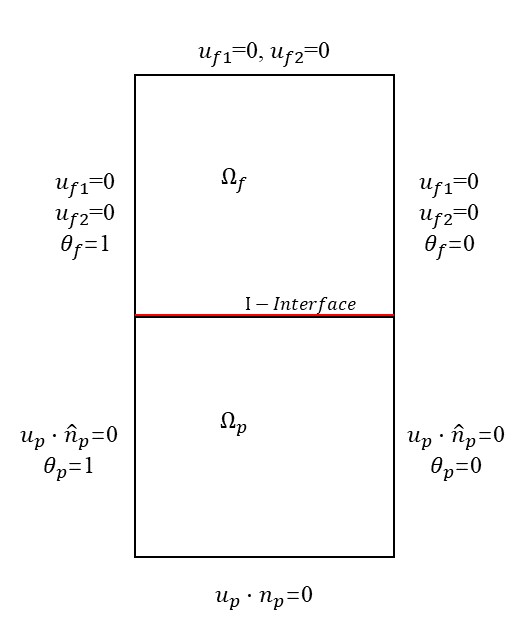}
  \end{minipage}%
}
  \centering
  \subfigure[]{
  \begin{minipage}[t]{0.3\linewidth}
  \centering
  \includegraphics[width=1.5in]{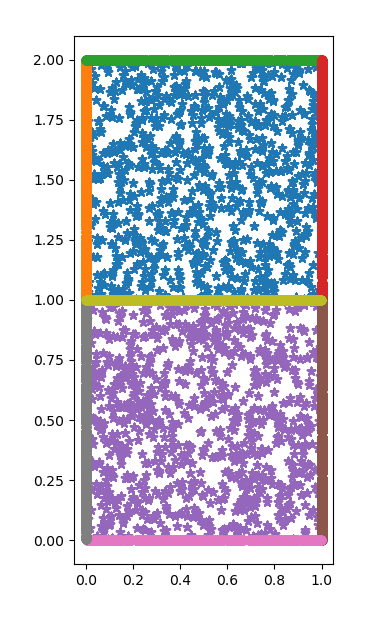}
  \end{minipage}%
}
  \caption{(a):~The computational domain with boundary conditions; (b):~The samples in the test domain.}\label{figure-figure6}
\end{figure}

Assume the computational domain is $\Omega_f=[0,1]\times[1,2]$, and $\Omega_p=[0,1]\times[0,1]$. The boundary conditions are illustrated in the Figure \ref{figure-figure6}. No-slip boundary condition for pipe region and no-flow boundary condition on the reservoir region are imposed for all the boundaries:
\begin{align}
u_f=0,~on~\partial\Omega_f\setminus\mathbb{I},~and~u_p\cdot n_p=0~on~\partial\Omega_p\setminus\mathbb{I}.\label{equation-4.2}
\end{align}
Left boundary is considered as heated wall with $\theta_f=\theta_p=1$ and right boundary is cold with $\theta_f=\theta_p=0$. On the interface between the two regions, the interface conditions \ref{equation-2.17}-\ref{equation-2.20}, which are proposed for the model, are utilized.

\begin{figure}[H]
  \centering
  \includegraphics[scale=0.4]{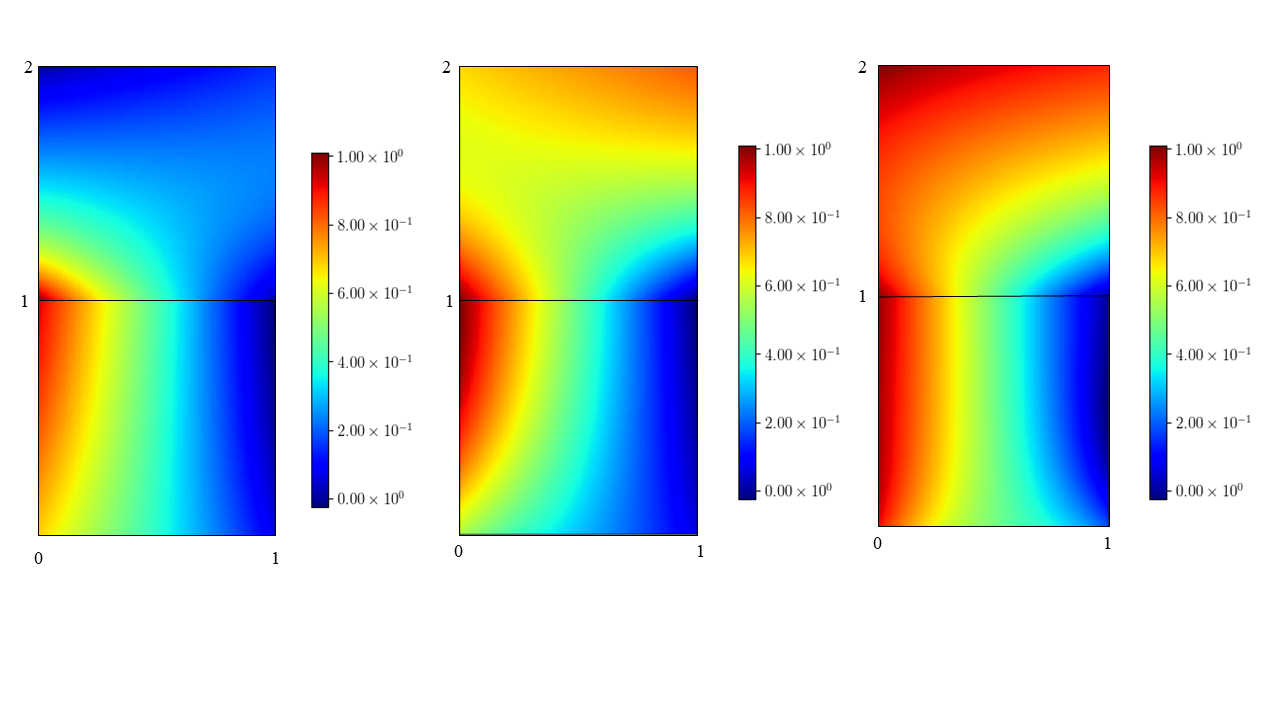}\\
  \caption{The temperature distribution for different values of Rayleigh number. Left: $Ra=10^3$; Middle: $Ra=10^4$; Right: $Ra=10^5$.}\label{figure-figure7}
\end{figure}

The parameters are chosen as $k_f=k_p=1.0$, $Pr=0.71$, $Ra=10^3\leq10^5$, and $Da=1.0\times10^{-2}$. The initial values for velocity and temperature in both subdomains are chosen to be $0$ and $1$ respectively. From Figure \ref{figure-figure7}, it is obvious that the thermal transfers from hot wall to cold wall. What's more, the thermal convection changes more violently as the Rayleigh number increases.


\end{document}